\newcommand{\R}{I\!\!R}   
\newcommand{\Z}{{\mathbb Z}} 

\newcommand{\sgn}{\mbox{\sl sgn}}
\newcommand{\dgn}{\mbox{\sl dgn}}
\newcommand{\rk}{\mbox{\sl rk}}
\newcommand{\mul}{\mbox{\sl mul}\,}
\newcommand{\Bsym}{{\mathrm B}_{\mathrm{sym}}}
\newcommand{\maslov}{\mathrm i_{\mathrm{maslov}}}
\newcommand{\iip}[2]{\langle#1,#2\rangle}
\newcommand{\op}{\mathrm{op}}

\newcommand{\sextuple}{(\mathcal M,\omega,H,\mathfrak L,\Gamma,\mathcal P)}
\newcommand{\Hf}{\vec H}
\newcommand{\HL}{H_{\mathfrak L}}

\newcommand{\Ddt}{\frac{\mathrm D}{\mathrm dt}}
\newcommand{\ddt}{\frac{\mathrm d}{\mathrm dt}}

\newcommand{\dds}{\frac{\mathrm d}{\mathrm ds}}
\newcommand{\ddp}{\frac{\partial}{\partial p}}
\newcommand{\ddq}{\frac{\partial}{\partial q}}
\newcommand{\vertical}{_{\mathrm{vert}}}

\documentclass[oneside,draft,letterpaper,10pt]{amsart}

\usepackage{amsmath}               
\usepackage{amsthm}                
\usepackage{times}

\numberwithin{equation}{section}

\title[An Index Theorem for Hamiltonian Systems]%
{An Index Theorem for non Periodic Solutions of Hamiltonian Systems}
\author[P.\ Piccione, D.\ Tausk]{Paolo Piccione\and Daniel V.\ Tausk}
\address{Departamento de Matem\'atica,   Universidade de S\~ao Paulo, Brazil}
\email{piccione@ime.usp.br, tausk@ime.usp.br}
\urladdr{http://www.ime.usp.br/\~{}piccione, http://www.ime.usp.br/\~{}tausk}
\thanks{The first author is partially sponsored by CNPq, the second author is
sponsored by FAPESP}

\date{revised version of June 2000}

\begin{document}


\theoremstyle{plain}\newtheorem{teo}{Theorem}[subsection]
\theoremstyle{plain}\newtheorem{prop}[teo]{Proposition}
\theoremstyle{plain}\newtheorem{lem}[teo]{Lemma} 
\theoremstyle{plain}\newtheorem{cor}[teo]{Corollary} 
\theoremstyle{definition}\newtheorem{defin}[teo]{Definition} 
\theoremstyle{remark}\newtheorem{rem}[teo]{Remark} 
\theoremstyle{plain} \newtheorem{assum}[teo]{Assumption}
\theoremstyle{definition}\newtheorem{example}[teo]{Example}


\begin{abstract}
We consider a {\em Hamiltonian setup\/} $\sextuple$,
where $(\mathcal M,\omega)$ is a symplectic manifold, 
$\mathfrak L$ is a distribution of Lagrangian subspaces in $\mathcal M$,
$\mathcal P$ a Lagrangian
submanifold of $ \mathcal M$, $H$ is a smooth time dependent
Hamiltonian function  on $\mathcal M$ and $\Gamma:[a,b]\to\mathcal M$ is an integral
curve  of the Hamiltonian flow $\Hf$ starting at $\mathcal P$. 
We do not require any convexity property of the Hamiltonian function $H$.
Under the assumption that $\Gamma(b)$ is not $\mathcal P$-focal it is introduced 
the Maslov index $\maslov(\Gamma)$ of $\Gamma$ given in terms of the first
relative homology group of the Lagrangian Grassmannian; under generic circumstances
$\maslov(\Gamma)$ is computed as a sort of {\em algebraic count\/}
of the $\mathcal P$-focal points along $\Gamma$. 
We prove the following version of the Index Theorem: under suitable hypotheses, the Morse
index of the Lagrangian action functional restricted to suitable variations
of $\Gamma$  is equal to the sum of $\maslov(\Gamma)$ and a {\em convexity term\/} of the
Hamiltonian
$H$ relative to the submanifold $\mathcal P$.  When the result is
applied to the case of the cotangent bundle $\mathcal M=TM^*$ of a semi-Riemannian manifold
$(M,g)$ and to the geodesic Hamiltonian $H(q,p)=\frac12 g^{-1}(p,p)$, we obtain
a semi-Riemannian version of the celebrated Morse Index Theorem for geodesics
with variable endpoints in Riemannian geometry.
\end{abstract}

\maketitle

\begin{section}{Introduction}\label{sec:intro}
Our interest in the index theory for solutions of the Hamilton 
equations in a symplectic manifold was originally
motivated by the aim of extending to the case of non positive definite
metrics the classical results of the Morse theory for geodesics in
Riemannian manifolds (see~\cite{M}). Despite this original motivation,
the geometric applications of the theory developed are left to the very
last part of the article, and most of the results presented in the paper
belong indeed to the realm of the theory of systems of ordinary differential
equations with coefficients in the Lie algebra 
$\mathrm{sp}(2n,\R)$ of the symplectic group.
Such systems will be called {\em symplectic differential systems}.

In order to motivate the theory presented in this paper,
we give a short account of the mathematical history 
of the problem. The origin of the index theory is to be found
in the {\em Sturm theory\/} for ordinary differential equations
(see for instance \cite{CL}).
The Sturm oscillation theorem deals
with second order differential equations of the form $-(px')'+rx=\lambda x$
where $p$ and $r$ are functions with $p>0$, and $\lambda$ is a real parameter.
%
The theorem states that the number of zeroes of a non null solution $x$ of
the Sturm equation satisfying $x(a)=0$ equals the index of the index form
$I(x_1,x_2)=\int_a^b[px_1'x_2'+r x_1x_2]\;{\rm d}t$ defined in
space of real valued maps on $[a,b]$ vanishing at the endpoints.

The extension of the results of the Sturmian theory to the case of systems
of differential equations is essentially due to Morse, obtaining the celebrated
Morse Index Theorem in Riemannian geometry (see for instance \cite{dC, M}). 
The Morse--Sturm systems for which the theorem applies are those of the form $g^{-1}(gv')'=Rv$,
where $g(t)$ is a {\em positive definite\/} symmetric matrix and $R(t)$ is $g(t)$-symmetric
linear operator on $\R^n$ for all $t$. Such systems are obtained, for instance, 
by considering the Jacobi equation along a geodesic in a Riemannian manifold;
the equation is converted into a system of ODE's in $\R^n$ by means of a parallel
trivialization of the tangent bundle of the manifold along the geodesic.
In this situation, the {\em index form\/} $I(v,w)=\int_a^b\big[g(v',w')+g(Rv,w)\big]\;\mathrm
dt$ has finite index in the space of vector valued function
on $[a,b]$ vanishing at the endpoints, and it is equal to the number of conjugate points of
the system in the interval $]a,b[$. 
By minor changes, the Morse Index Theorem is also valid in the case
of {\em Lorentzian\/} metrics $g$, i.e., metrics having index one, 
provided that one considers causal, i.e., timelike or lightlike, geodesics 
and that one restricts the index form $I$ to vector fields that are
pointwise orthogonal to the geodesic (see \cite{BEE}).
Subsequent results have extended the Index Theorem to the case
of solutions with non fixed initial and/or final endpoint (see
\cite{Kal, PT})  and to the case of ordinary differential operators
of even order (see~\cite{Ed}).

When passing to the case of spacelike geodesics in Lorentzian manifolds,
or, more in general, to geodesics in semi-Riemannian manifolds endowed
with metric tensors of arbitrary index, there is no hope to extend the
original formulation of the Morse Index Theorem for several reasons.
In first place,  the action functional is strongly indefinite, i.e., 
its second variation, given
by the index form $I$, has {\em always\/} infinite index. Moreover, the
set of conjugate points along a given geodesic may fail to be discrete,
and the Jacobi differential operator is no longer self-adjoint.

A different integer valued geometric invariant, called the {\em Maslov index}, has been recently
introduced in the context of semi-Riemannian
geodesics (see \cite{GMPT, Hel1, MPT}). 

In the case of Riemannian or causal
Lorentzian geodesics, the Maslov index coincides with the {\em geometric
index\/} of the geodesic, which is the number of conjugate points counted
with multiplicity. For metrics of arbitrary index, under generic
circumstances, this index is computed as a sort of algebraic count of the conjugate points
along the geodesic.  
It is natural to expect that, as in the case of positive definite metrics
(see \cite{Duis}), the Maslov index should
play the role of the geometric index in metrics with arbitrary index, and in this paper we
present several arguments to strengthen this idea. Namely, we show that
the Maslov index, under certain circumstances,    is  equal to the index
of the restriction of $I$ to a suitable subspace of variational
vector fields. This equality gives a generalization of the Morse
Index Theorem to general indefinite metrics.

The definition of the Maslov index is obtained by developing
an intersection theory for curves in the manifold $\Lambda$ of all
Lagrangians of a symplectic space. 
A non trivial solution of the symplectic differential system gives a continuous curve
in the Lagrangian Grassmannian, and the zeroes of the solution correspond
to intersections of this curve with the subvariety $\Lambda_{\ge1}(L_0)$
of all Lagrangians
$L$ which are not complementary to a fixed Lagrangian $L_0$. The Maslov index $\maslov$
of the symplectic differential system is then defined using the first relative
homology group of the pair $(\Lambda,\Lambda\setminus\Lambda_{\ge1}(L_0))$.
There is no standard notion of Maslov index in the literature, which is defined differently
according to the context. We use the definition given in \cite{Ar, Duis, Hel1, MPT};
we remark that for periodic Hamiltonian systems a different
definition is usually adopted (see \cite{Ek, MW, SZ}).

Symplectic systems arise naturally as linearizations of the
Hamilton equations. We consider a general Hamiltonian setup, consisting
of a symplectic manifold $(\mathcal M,\omega)$, a Lagrangian submanifold $\mathcal P$
of $\mathcal M$, a distribution $\mathfrak L$ of Lagrangian spaces in $\mathcal M$,
a time dependent Hamiltonian function
$H$ on $\mathcal M$ and a given solution $\Gamma:[a,b]\to\mathcal M$
of the Hamilton  equations, with $\Gamma(a)\in\mathcal P$. 
We introduce an index form $I_\Gamma$ associated to these objects;
for instance, if we consider the cotangent bundle $\mathcal M=TM^*$
of a manifold $M$, where the distribution
$\mathfrak L$ is the vertical bundle of $TTM^*$, 
and $H$ is a hyper-regular Hamiltonian, then our index form $I_\Gamma$
coincides with the second variation of the Lagrangian action functional.
If $M$ is endowed with a semi-Riemannian metric $g$ and $H$ is given
by $H(q,p)=\frac12g^{-1}(p,p)$, then the solution $\Gamma$ projects
onto a geodesic $\gamma$ in $(M,g)$, and $I_\Gamma$ becomes the standard
index form of semi-Riemannian geodesics.

A genuine physical  interest in the Hamiltonian extension
of the index theory comes from Mechanics and Optics, both
classical and relativistic, where the Hamiltonian formalism
appears naturally in many situations. For instance, the authors 
are currently studying a Hamiltonian 
formalism for light rays in a relativistic medium (see \cite{PPT},
see also \cite{Perlick1999} for a complete description of the
relativistic ray optics). As the underlying spacetime 
model it is assumed an arbitrary 4-dimensional Lorentzian manifold, 
and the medium is described by a non convex Hamiltonian function  
that typically involves the spacetime metric and 
a number of tensor fields by which the
medium is characterized. 

Given a semi-Riemannian geodesic $\gamma:[a,b]\to M$,
then by a parallel trivialization of the tangent bundle $TM$ along $\gamma$
the Jacobi equation becomes a Morse--Sturm system. Similarly,
given a Hamiltonian setup $\sextuple$, then by using a symplectic referential
of $T\mathcal M$ along $\Gamma$ adapted to the distribution $\mathfrak L$,
we obtain a symplectic differential system, that corresponds to the
linearization of the Hamilton equations along $\Gamma$.
Our index form $I_\Gamma$ will then correspond to a symmetric bilinear form $I$
in the space of maps $v:[a,b]\to\R^n$ satisfying suitable boundary
conditions. The main result of this paper is an index theorem (Theorem~\ref{thm:indextheorem})
that relates the index of a suitable restriction of $I$ (or, equivalently, of $I_\Gamma$) with the
Maslov index of a symplectic differential system.
Such results aims at the developments of an infinite dimensional Morse
theory for solutions of Hamiltonian systems, in the spirit of
\cite{GMPT}, where the authors obtain existence results for geodesics
in stationary Lorentzian manifolds.

The Hessian of the Hamiltonian function $H$ with respect to the momenta
corresponds to a component $B$ of the coefficient matrix of the symplectic
differential system.
In the case that the Hamiltonian is convex, i.e., when $B$ is positive definite, 
such index theorem is essentially the classical Morse Index Theorem for Riemannian geodesics,
or equivalently, for Morse--Sturm systems with $g$ positive definite.
In the non convex case, the index of $I$ is always infinite; in order to obtain
a space $\mathcal K$ carrying a finite index, we ``factor out'' an infinite
dimensional space $\mathcal S$ from the space of all variations, such that
$I$ is negative definite on $\mathcal S$. The spaces $\mathcal K$
and $\mathcal S$ are orthogonal with respect to the index form $I$.
The definition of the spaces $\mathcal K$ and $\mathcal S$ is based
on the choice of a {\em maximal negative distribution\/} for $B$, i.e.,
a family $\{\mathcal D_t\}_{t\in[a,b]}$ of subspaces of $\R^n$ such
that each $\mathcal D_t$ is a maximal negative subspace for $B(t)$.
The space $\mathcal K$ is essentially the space of variations $v:[a,b]\to\R^n$
that are ``solutions of the symplectic differential system  in the directions 
of $\mathcal D$'', while the space $\mathcal S$ consists of variations taking values
in $\mathcal D$. A geometric description of this abstract setup in the 
context of semi-Riemannian geodesics is given in the examples at the
end of Section~\ref{sec:hamiltonian}.
The intersection $\mathcal K\cap\mathcal S$ consists of solutions
of another symplectic differential system, called the {\em reduced
symplectic system}, that vanish at the endpoints.
The main technical hypotheses for our theory  are a nondegeneracy assumption
on the initial condition and the assumption that the reduced symplectic
system has no conjugate points on the entire interval $]a,b]$.
An even more general result could be proven
if this latter condition is dropped, provided that one takes into account
in the thesis of the index theorem also the contribution
given by the conjugate points of the reduced symplectic system.
Some results in this direction were announced in \cite{CRASP}.

There exist in the literature other   index theories for  
Hamiltonian systems, that mostly concern the case of periodic
solutions. 
Theorems (like the one presented in this paper) relating a version of the Maslov index of the
solution with the index of the second variation of the Lagrangian action functional
were proven only in the case of convex Hamiltonians (see for instance \cite{Ek, MW}).
In \cite{CZ} and \cite{Long} the authors develop a Morse theory for  periodic solutions of
non convex Hamiltonian systems employing a technique that reduces the problem
to finite dimensional Morse theory.   The index theorems presented in \cite{Long, CZ}
relate the Maslov index of a periodic solution to the Hessian
of suitable finite-dimensional versions of the action functional.

The paper is divided into two sections. In Section~\ref{sec:symplectic} we
study the theory of symplectic differential systems in $\R^n$ and in
Section~\ref{sec:hamiltonian} we discuss the applications to Hamiltonian
systems and semi-Riemannian geometry.
In order to facilitate the reading of the article, we give a short overview
of the contents of each subsection.
In subsection~\ref{sub:basic} we introduce the notion of symplectic differential
system, with the appropriate initial conditions and we define the basic
notions, like focal instants and focal index. In subsection~\ref{sub:second}
we state a criterion that says when a second-order linear differential
equation in $\R^n$ can be transformed into a symplectic differential system.

In subsection~\ref{sub:maslov} we describe briefly the geometry of the
Lagrangian Grassmannian of a symplectic space and we define the {\em Maslov
index\/} associated to a symplectic differential problem. We show that the
Maslov index equals the focal index, under a certain non-degeneracy
assumption. In subsection~\ref{sub:nondegeneracy} we show that this
non-degeneracy condition is {\em generic}, i.e., every symplectic problem can
be uniformly approximated by symplectic problems satisfying the condition.
In subsection~\ref{sub:indexform} we define the index form $I$ associated to a
symplectic differential problem and we show that, except for the convex case, 
this index form has
always infinite index. In subsection~\ref{sub:K} we describe the space $\mathcal K$  and we
introduce the reduced symplectic system. In subsection~\ref{sub:indextheorem} we compute the index
of $I$ in $\mathcal K$, proving the index theorem.
In subsection~\ref{sub:variable}, motivated by the index theorem for geodesics
with both endpoints variable in submanifolds (see \cite{Kal,PT}), we extend
the index theorem of subsection~\ref{sub:indextheorem} to the case of
symplectic differential problems with boundary conditions at both endpoints of
the interval.
In subsection~\ref{sub:opposite} we  obtain an
alternative version of the index theorem which involves the index of $-I$,
by introducing the notion of {\em opposite\/} symplectic system.
In subsection~\ref{sub:equivalence} we define a notion of isomorphism for
symplectic differential systems, motivated by the fact that Hamiltonian
systems determine symplectic differential systems only up to isomorphism,
depending on the choice of the symplectic trivialization along the
solution. We prove that all the objects of our theory
(focal instants, Maslov index, index form) are invariant by such isomorphisms. In
subsection~\ref{sub:isoMS} we show that every symplectic differential system
is isomorphic to a Morse--Sturm system.
In subsection~\ref{sub:associated} we describe how to produce a symplectic differential
problem from a Hamiltonian setup by using a symplectic trivialization
of the tangent bundle of the symplectic manifold along a solution.
It is then possible, thanks to the results of subsection~\ref{sub:equivalence},
to define the Maslov index and an index form
for a Hamiltonian setup. In subsection~\ref{sub:maslovH} we restate the index
theorem in terms of the Hamiltonian setup.
In subsection~\ref{sub:hyper} we show that, in the case of a hyper-regular
Hamiltonian, the index form coincides with the second variation of the action
functional associated to the corresponding Lagrangian. 
In subsection~\ref{sub:interpretation} we describe the Hamilton 
equations and their linearized version in terms of a connection in the base manifold $M$.
Finally, in subsection~\ref{sub:semiRiemannian} we give a geometric application
of the index theorem to semi-Riemannian geometry.

\end{section}


\begin{section}{Symplectic Differential Systems}
\label{sec:symplectic}
Throughout the paper we use the following convention: given
finite dimensional real vector spaces $V$ and $W$, we will denote
by $W^*$ the dual space of $W$ and, given a linear
map $T:V\to W$, we will identify it with the bilinear map
$T:V\times W^*\to\R$ given by $T(v,\alpha)=\alpha(T(v))$.
From the context, there should be no confusion arising from using the same symbol
for a linear map and the corresponding bilinear map.
The reader should be warned that, for a fixed choice of basis in $V$ and $W$,
the matrix representations of a linear map $T$ and of the corresponding bilinear map
do {\em not\/} coincide, but they are the transpose of each other.
If $V=W$, we observe that the linear map $T$ is invertible if and only if
the corresponding bilinear map is nondegenerate.
We denote by $T^*:W^*\to V^*$  the adjoint map $T^*(\alpha)=\alpha\circ T$;
observe that, in terms of bilinear maps, one has $T(v,\alpha)=T^*(\alpha,v)$
for all $v\in V$ and $\alpha\in W^*$.

As a consequence of these identifications, we have that, if $B:V\times W\to\R$
is a bilinear map, and $T:Z\to V$, $S:Z\to W$ are linear maps on some
vector space $Z$, then the composite linear maps $B\circ T$ and $S^*\circ B$
correspond respectively to the bilinear maps $B(T\,\cdot\,,\,\cdot\,)$ and
$B(\,\cdot\,,S\,\cdot\,)$.

Given Banach spaces $E_1$ and $E_2$, we denote by $\mathcal L(E_1,E_2)$
the set of all bounded linear operators from $E_1$ to $E_2$ and by $\mathrm B(E_1,E_2,\R)$
the set of all bounded bilinear maps from $E_1\times E_2$ to $\R$. If $E_1=E_2=E$, we also
set $\mathcal L(E)=\mathcal L(E,E)$ and $\mathrm B(E,\R)=\mathrm B(E,E,\R)$; by $\mathrm
B_{\mathrm{sym}}(E,\R)$ we mean the set of symmetric bounded bilinear maps on $E$.
If $E=H$ is a Hilbert space with a fixed inner product $\iip\cdot\cdot$, by Riesz's
theorem we can identify $\mathcal L(H)$ with $\mathrm B(H,\R)$ by associating
the linear operator $T$ to the bilinear form $\iip{T\cdot}{\cdot}$; obviously,
the bilinear form is symmetric if and only if the linear operator is self-adjoint.

We give some general definitions 
concerning symmetric bilinear forms for later use.
 Let $V$ be any real vector space and $B:V\times V\to\R$ a symmetric
bilinear form. The {\em negative type number\/} (or {\em index}) $n_-(B)$ of $B$
is the possibly infinite number defined by
\begin{equation}\label{eq:def-}
n_-(B)=\sup\big\{{\rm dim}(W):W\ \text{subspace of}\ V \ \text{on which}\ B\ 
\text{is negative definite}\big\}.
\end{equation}
The {\em positive type number\/} (or {\em co-index}) $n_+(B)$ is given by
$n_+(B)=n_-(-B)$; if at least one of these two numbers is finite,
the {\em signature\/} $\sgn(B)$ is defined by:
\[\sgn(B)=n_+(B)-n_-(B).\]
The {\em kernel\/} of $B$, ${\rm Ker}(B)$, is the set  of vectors $v\in V$
such that $B(v,w)=0$ for all $w\in V$; the {\em degeneracy\/} $\dgn(B)$
of $B$ is the (possibly infinite) dimension of ${\rm Ker}(B)$. 
If $V$ is finite dimensional, then the numbers $n_+(B)$, $n_-(B)$ and
$\dgn(B)$ are respectively  the number of $1$'s, $-1$'s and $0$'s
in the canonical form of $B$ as given by the Sylvester's Inertia Theorem.
In this case, $n_+(B)+n_-(B)$ is equal to the codimension of ${\rm Ker}(B)$,
and it is also called the {\em rank} of $B$, $\rk(B)$.

\subsection{Basic definitions}\label{sub:basic}
We consider the symplectic space $\R^n\oplus{\R^n}^*$ endowed
with the canonical symplectic form $\omega$:
\begin{equation}\label{eq:defomega}
\omega((v_1,\alpha_1),(v_2,\alpha_2))=\alpha_2(v_1)-\alpha_1(v_2).
\end{equation}
We denote by $\mathrm{Sp}(2n,\R)$ the Lie group of symplectic transformations
of the space $(\R^n\oplus{\R^n}^*,\omega)$ and by $\mathrm{sp}(2n,\R)$ its Lie algebra.
Recall that an element $X\in\mathrm{sp}(2n,\R)$ is a linear map
$X:\R^n\oplus{\R^n}^*\to\R^n\oplus{\R^n}^*$ such that $\omega(X\,\cdot\,,\,\cdot\,)$ is symmetric;
in block matrix form, $X$ is given by:
\begin{equation}\label{eq:X}
X=\left(\begin{array}{cc}A&B\\ C&-A^*\end{array}\right),
\end{equation}
where $A:\R^n\to\R^n$ is an arbitrary linear map, and $B:{\R^n}^*\to\R^n$,
$C:{\R^n}\to{\R^n}^*$ are symmetric when regarded as bilinear maps.

The main object of our study are the following differential
systems:
\begin{defin}\label{thm:defsympl}
A {\em symplectic differential system\/} in $\R^n$ is a first 
order linear differential system of the form:
\begin{equation}\label{eq:sympl}
\left\{\begin{array}{l}v'(t)=A(t)v(t)+B(t)\alpha(t);\\\alpha'(t)=C(t)v(t)-A^*(t)\alpha(t),
\end{array}\right.\quad t\in [a,b],\ v(t)\in\R^n,\ \alpha(t)\in{\R^n}^*
\end{equation}
where the coefficient matrix \[X(t)=\left(\begin{array}{cc}A(t)&B(t)\\
C(t)&-A^*(t)\end{array}\right)\]
is a  curve in $\mathrm{sp}(2n,\R)$ (i.e., $B(t)$ and $C(t)$ are symmetric), 
with $A$ and $B$ of class $C^1$,
$C$ continuous  and $B(t)$ invertible for all $t$. Some of the results
will be proven under the more restrictive assumption that $C$ be a map
of class $C^1$ (see Subsection~\ref{sub:indextheorem}), however, for the
final results, we will get rid of this extra hypothesis by a perturbation
argument. Some results concerning the isomorphisms of symplectic systems
will be proven under the assumption that $B$ is a map of class $C^2$
(Subsection~\ref{sub:isoMS}).
\end{defin}
Given a symplectic differential system \eqref{eq:sympl}, we denote
by $\Psi(t)$, $t\in[a,b]$, its {\em fundamental matrix}, defined by:
\begin{equation}\label{eq:defPsi}
\Psi(t)(v(a),\alpha(a))=(v(t),\alpha(t)),
\end{equation}
for all solution $(v(t),\alpha(t))$ of \eqref{eq:sympl}. The map $\Psi$
can be characterized as the curve in the Lie group of linear isomorphisms
of $\R^n\oplus{\R^n}^*$ satisfying:
\begin{equation}\label{eq:charPsi}
\Psi'(t)=X(t)\,\Psi(t),\quad \Psi(a)=\mathrm{Id}.
\end{equation}
Since $X(t)\in\mathrm{sp}(2n,\R)$, then $\Psi(t)\in\mathrm{Sp}(2n,\R)$
for all $t$; this means that for any pair of solutions $(v,\alpha)$ and $(w,\beta)$
of \eqref{eq:sympl} we have:
\begin{equation}\label{eq:abvwcosnt}
\beta(v)-\alpha(w)=\text{constant}.
\end{equation}
Given any $C^1$-curve $v:[a,b]\to\R^n$, the first equation
of \eqref{eq:sympl} defines a unique continuous
curve $\alpha_v:[a,b]\to{\R^n}^*$;
more explicitly:
\begin{equation}\label{eq:defalphav}
\alpha_v(t)=B(t)^{-1}\Big(v'(t)-A(t)v(t)\Big).
\end{equation}
For future reference, we remark the following equality:
\begin{equation}\label{eq:alphafv}
\alpha_{fv}=f'B^{-1}(v)+f\alpha_v,
\end{equation}
for any $C^1$-function $f:[a,b]\to\R$.

Let $\ell_0$ be {\em Lagrangian subspace\/} of $(\R^n\oplus{\R^n}^*,\omega)$.
This means that $\mathrm{dim}(\ell_0)=n$ and $\ell_0$ is {\em $\omega$-isotropic}, i.e., 
$\omega$ vanishes on $\ell_0\times\ell_0$. 
We consider the following initial conditions for the system
\eqref{eq:sympl}:
\begin{equation}\label{eq:IC}
(v(a),\alpha(a))\in\ell_0.
\end{equation}
There exists a bijection between the set of Lagrangian subspaces
of $(\R^n\oplus{\R^n}^*,\omega)$ and the set of pairs $(P,S)$, where
$P$ is a subspace of $\R^n$ and $S:P\times P\to\R$ is a symmetric bilinear
form. The bijection is defined by:
\begin{equation}\label{eq:bijection}
\ell_0=\big\{(v,\alpha):v\in P,\ \alpha\vert_P+S(v)=0\big\}.
\end{equation}
We can therefore rewrite the initial conditions \eqref{eq:IC} in terms of $(P,S)$:
\begin{equation}\label{eq:ICPS}
v(a)\in P,\quad\alpha(a)\vert_P+S(v(a))=0.
\end{equation}
For the theory developed in this paper we will make henceforth the following:
\begin{assum}\label{thm:nondeg}
We assume that the bilinear form $B(a)^{-1}$ in $\R^n$
is nondegenerate on $P$, or, equivalently, that $B(a)$ is non degenerate
on the annihilator $P^o$ of $P$, which is easily seen to be given by:
\[P^o=\big\{\alpha\in{\R^n}^*:(0,\alpha)\in \ell_0\big\}.\]
\end{assum}

A solution $(v,\alpha)$ of \eqref{eq:sympl} will be called an {\em $X$-solution\/};
if $(v,\alpha)$ in addition satisfies \eqref{eq:IC}, it will
be called an $(X,\ell_0)$-solution. A pair $(X,\ell_0)$ as described above
satisfying Assumption~\ref{thm:nondeg} will be called a {\em set of data for
the symplectic differential problem}. 
By a slight abuse of terminology, we will say that a $C^2$-curve
$v(t)$ is an $X$-solution if the pair
$(v,\alpha_v)$, with $\alpha_v$ given by \eqref{eq:defalphav}, 
is an  $X$-solution; similarly, we say that  
$v$ is an $(X,\ell_0)$-solution  
if the pair $(v,\alpha_v)$ is an $(X,\ell_0)$-solution.

We now consider a fixed set of data $(X,\ell_0)$ for the symplectic differential
problem. We denote by $\mathbb V$ the set of all $(X,\ell_0)$-solutions:
\begin{equation}\label{eq:spaces}
\mathbb V=\big\{v: v\ \text{is an $(X,\ell_0)$-solution}\big\}.
\end{equation}
The fact that $\ell_0$ is Lagrangian, combined with \eqref{eq:abvwcosnt}
implies that 
\begin{equation}\label{eq:abvwcosnt2}
\alpha_v(w)=\alpha_w(v),\quad\forall\,v,w\in\mathbb V.
\end{equation}
For $t\in[a,b]$, we set 
\[\mathbb V[t]=\big\{v(t):v\in\mathbb V\big\}.\]
From \eqref{eq:abvwcosnt2} and a simple dimension counting argument,
the annihilator of $\mathbb V[t]$ in ${\R^n}^*$ is given by:
\begin{equation}\label{eq:VTo}
\mathbb V[t]^o=\big\{\alpha_v(t):v\in\mathbb V,\ v(t)=0\big\},\quad t\in\,[a,b]. 
\end{equation} 
\begin{defin}\label{thm:deffocal}
An instant $t\in\,]a,b]$ is said to be {\em focal\/} for the pair
$(X,\ell_0)$ if there exists a non zero $v\in\mathbb V$ 
such that $v(t)=0$, i.e., if $\mathbb V[t]\ne\R^n$.
The {\em multiplicity\/} $\mul(t)$ of the focal instant $t$ is defined to
be the dimension of the space of those $v\in\mathbb V$ vanishing at $t$, or,
equivalently, the codimension of $\mathbb V[t]$ in $\R^n$.
The {\em signature} $\sgn(t)$ of the focal instant $t$ is the signature
of the restriction of the bilinear form $B(t)$ to the space $\mathbb V[t]^o$, or,
equivalently, the signature of the restriction of $B(t)^{-1}$ to
the $B(t)^{-1}$-orthogonal complement $\mathbb V[t]^\perp$ of $\mathbb V[t]$ 
in $\R^n$.
The focal instant $t$ is said to be {\em nondegenerate\/} if such restriction
is nondegenerate.
If there is only a finite number of focal instants in $]\,a,b]$, we define
the {\em focal index\/} $\mathrm i_{\mathrm{foc}}=\mathrm i_{\mathrm{foc}}(X,\ell_0)$ 
to be the sum:
\begin{equation}\label{eq:defifocal}
\mathrm i_{\mathrm{foc}}=\sum_{t\in\,]a,b]}\sgn(t).
\end{equation}
\end{defin}
\begin{rem}\label{thm:remrealanalytic}
There exist several situations where the number of focal instants is indeed
finite; for instance, this is always the case when  the coefficient matrix
$X$ is real analytic in $t$. This fact follows from the observation
that the focal instants are the zeroes of the function $\mathrm{det}(v_1(t),\ldots,v_n(t))$,
where $\{v_1,\ldots,v_n\}$ is a basis of $\mathbb V$. Such function does not vanish
identically, as we will see in Subsection~\ref{sub:maslov} that it does not vanish in a
neighborhood of $t=a$. We will prove 
in  Subsection~\ref{sub:maslov} that nondegenerate focal instants are isolated.
\end{rem}

In a sense, Assumption~\ref{thm:nondeg} says that $t=a$ is a nondegenerate focal
instant; observe that $\mathbb V[a]=P$.
\subsection{Morse--Sturm and second order linear differential equations}
\label{sub:second}
An important class of examples of symplectic differential systems
arises from the so called {\em Morse--Sturm systems}, which are second
order differential equations in $\R^n$ of the form:
\begin{equation}\label{eq:MS}
g^{-1}(gv')'=Rv,
\end{equation}
where $g(t)$ is a $C^1$-curve of symmetric nondegenerate bilinear forms
in $\R^n$ and $R(t)$ is a continuous curve of linear operators
on $\R^n$ such that $g(t)(R(t)\cdot,\cdot)$ is symmetric for all $t$.
Considering the change of variable $\alpha=gv'$, the equation
\eqref{eq:MS} is equivalent to the following symplectic system:
\begin{equation}\label{eq:MSeq}
\left\{\begin{array}{l}v'=g^{-1}\alpha,\\ \alpha'=gRv. \end{array}\right.
\end{equation}
Comparing with \eqref{eq:sympl}, in \eqref{eq:MSeq} we have $A=0$, $B=g^{-1}$, $C=gR$.

In general, \eqref{eq:sympl} may be written as a second order equation, as follows.
Substitution of \eqref{eq:defalphav} into the second equation of
\eqref{eq:sympl} shows that $v$ is a solution of \eqref{eq:sympl}
if and only if it is a solution of the following  second order equation:
\begin{equation}\label{eq:secord}
\left[B^{-1}(v'-Av)\right]'=Cv-A^*B^{-1}(v'-Av).
\end{equation}
 One has the following natural question: when is the second order linear equation
\begin{equation}\label{eq:question}
v''+Z_1v'+Z_2v=0,
\end{equation}
with $Z_1,Z_2:[a,b]\to\mathcal L(\R^n)$ continuous, of the form \eqref{eq:secord}
for some symplectic differential system? 
It is not hard to prove the following:
\begin{prop}\label{thm:aimed}
The differential equation \eqref{eq:question} arises from 
a symplectic system if and only if there exists a (fixed)
symplectic form $\overline\omega$ in $\R^{2n}$ such that
$\Omega(t)^{-1}(\{0\}\oplus\R^n)$ is Lagrangian in $(\R^{2n},\overline\omega)$
for all $t$, where $\Omega(t):\R^{2n}\to\R^{2n}$ is the fundamental
matrix of \eqref{eq:question}. 
\end{prop}
\begin{proof}
Assume that $\xi(t)=\Omega(t)^{-1}(\{0\}\oplus\R^n)$ is Lagrangian in $(\R^{2n},\overline\omega)$
for all $t$ and choose a $C^1$-family of symplectomorphisms 
$\phi(t):(\R^{2n},\overline\omega)\to(\R^n\oplus{\R^n}^*,\omega)$ such that
$\phi(t)(\xi(t))=\{0\}\oplus{\R^n}^*$ for all $t$. Define $\Psi(t)=\phi(t)\phi(a)^{-1}$ and
$X(t)=\Psi'(t)\Psi(t)^{-1}$; then $X$ is the coefficient matrix of a symplectic system that
corresponds to \eqref{eq:question} by the change of variables $(v(t),\alpha(t))=
\phi(t)\Omega(t)^{-1}(v(t),v'(t))$.
\end{proof}
\subsection{The Maslov Index}\label{sub:maslov}
The goal of this subsection is to produce an integer-valued invariant
for a set of data $(X,\ell_0)$ for the symplectic differential
problem. We call this number the {\em Maslov index\/}
of the pair $(X,\ell_0)$, we show that it is stable by $C^0$-small
perturbations of the data, and that it is {\em generically\/} equal
to the focal index $\mathrm i_{\mathrm{foc}}(X,\ell_0)$.
We will use several well known facts about the geometry  of the Lagrangian
Grassmannian of a symplectic space (see for instance \cite{Ar, Duis, GS, MPT, Tre}); 
in particular, we will make
full use of the notations and of the results proven in Reference~\cite{MPT}.

Recalling the definition of $\Psi$ in formula~\eqref{eq:defPsi},
we start by observing that the space:
\begin{equation}\label{eq:defellt}
\ell(t)=\Psi(t)(\ell_0)=\big\{(v(t),\alpha_v(t)):v\in\mathbb V\big\}
\end{equation}
is Lagrangian for all $t$. Let $\Lambda$ be the set of all Lagrangian
subspaces of $(\R^n\oplus{\R^n}^*,\omega)$; $\Lambda$ is a compact, connected
real analytic embedded $\frac12n(n+1)$-dimensional 
submanifold of the Grassmannian of all $n$-dimensional
subspaces of $\R^n\oplus{\R^n}^*$.  
For $L\in\Lambda$, the tangent space $T_L\Lambda$ can be canonically
identified with the space $\Bsym(L,\R)$ of all symmetric bilinear
forms on $L$.

Let $L_0\in\Lambda$
be fixed; we define the following subsets of $\Lambda$:
\begin{equation}\label{eq:defLambdak}
\Lambda_k(L_0)=\big\{L\in\Lambda:\mathrm{dim}(L\cap L_0)=k\big\},\quad k=0,\ldots,n.
\end{equation}
Each $\Lambda_k(L_0)$ is a connected embedded analytic submanifold of $\Lambda$
having codimension $\frac12k(k+1)$ in $\Lambda$; $\Lambda_0(L_0)$ is a dense
open contractible subset of $\Lambda$, while its complementary set:
\begin{equation}\label{eq:Lambda>=1}
\Lambda_{\ge1}(L_0)=\bigcup_{k=1}^n\Lambda_k(L_0)
\end{equation}
is {\em not\/} a regular submanifold of $\Lambda$. It is an algebraic variety,
and its regular part is
given by $\Lambda_1(L_0)$, which is a dense open subset of $\Lambda_{\ge1}(L_0)$.
Observe that $\Lambda_1(L_0)$ has codimension $1$ in $\Lambda$; moreover, it
has a {\em natural\/} transverse orientation in $\Lambda$, which is
canonically associated to the symplectic form $\omega$  
(see~\cite[Section~3]{MPT} for the proofs and the details of these results). 

The first relative singular homology group with integer
coefficients $H_1(\Lambda,\Lambda_0(L_0))$ is computed in \cite[Section~4]{MPT}:
\begin{equation}\label{eq:isomorphism}
H_1(\Lambda,\Lambda_0(L_0))\simeq\Z,
\end{equation}
where the choice of the above isomorphism is related to the choice
of a transverse orientation of $\Lambda_1(L_0)$ in $\Lambda$, and it
is therefore canonical. Every continuous curve $l$ in $\Lambda$ with
endpoints in $\Lambda_0(L_0)$ defines an element in $H_1(\Lambda,\Lambda_0(L_0))$,
and we denote by 
\[\mu_{L_0}(l)\in\Z\]
the integer number corresponding to the homology class of $l$ by
the isomorphism \eqref{eq:isomorphism}. This number, which is additive
by concatenation and invariant by homotopies with endpoints in
$\Lambda_0(L_0)$,  is
to be interpreted as an {\em intersection number\/} of the curve $l$ with
$\Lambda_{\ge1}(L_0)$ (see Proposition~\ref{thm:compmaslov} below).

Recalling formula \eqref{eq:defellt}, 
we note that $\ell$ is a $C^1$-curve in $\Lambda$;
we fix 
\[L_0=\{0\}\oplus{\R^n}^*\]
and we observe that the focal instants of the pair $(X,\ell_0)$ coincide
with the intersections of the curve $\ell$ with $\Lambda_{\ge1}(L_0)$.
Moreover, $\ell(t)\in\Lambda_k(L_0)$ if and only if $\mul(t)=k$.

We want to define the {\em Maslov index\/} $\maslov(X,\ell_0)$ to be
the integer $\mu_{L_0}(\ell)$; however, we observe that the curve $\ell$ may fail
to have endpoints in $\Lambda_0(L_0)$. In order to have $\ell(b)\in\Lambda_0(L_0)$,
we need to assume that $t=b$ is not focal. The initial point
$\ell(a)$ is not in $\Lambda_0(L_0)$ unless $P=\R^n$; however, we will prove
below that, because of  Assumption~\ref{thm:nondeg}, the curve
$\ell$ has at the most an isolated intersection with $\Lambda_{\ge1}(L_0)$
at $t=a$. We can therefore give the following:
\begin{defin}\label{thm:defmaslov}
Assume that $t=b$ is not a $(X,\ell_0)$-focal instant.
The Maslov index $\maslov=\maslov(X,\ell_0)$ is defined as:
\begin{equation}\label{eq:defmaslov}
\maslov=\mu_{L_0}(\ell\vert_{[a+\varepsilon,b]}),
\end{equation}
where $\varepsilon>0$ is small enough so that $\ell$ does not intercept
$\Lambda_{\ge1}(L_0)$ in $]a,a+\varepsilon]$.
\end{defin}
Clearly, the quantity on the right hand side of equality
\eqref{eq:defmaslov} does not depend on the choice of $\varepsilon$.

We now describe a method for computing $\mu_{L_0}(l)$ of a $C^1$-curve
$l$ in $\Lambda$ under some generic conditions. Recall that tangent
vectors to $\Lambda$ are canonically identified with symmetric bilinear
forms.
\begin{prop}\label{thm:compmaslov}
Let $l:[a,b]\to\Lambda$ be a curve of class $C^1$, and let $t_0\in[a,b]$
be such that $l(t_0)\in\Lambda_{\ge1}(L_0)$ and $l'(t_0)$ is non degenerate
on $l(t_0)\cap L_0$. Then, $t_0$ is an isolated intersection of $l$ with
$\Lambda_{\ge1}(L_0)$. Moreover, if $l$ has endpoints in $\Lambda_0(L_0)$
and all the intersections of $l$ with $\Lambda_{\ge1}(L_0)$ satisfy the above 
nondegeneracy condition, then $\mu_{L_0}(l)$ can be computed as:
\begin{equation}\label{eq:util}
\mu_{L_0}(l)=\sum_{t\in\,]a,b[}\sgn(l'(t)\vert_{l(t)\cap
L_0}).
\end{equation}
\end{prop}
\begin{proof}
See \cite[Corollary~4.3.3]{MPT}.
\end{proof}
Observe in particular that, under the assumptions of Proposition~\ref{thm:compmaslov},
if all the intersections of $l$ with $\Lambda_{\ge1}(L_0)$ lie in
$\Lambda_1(L_0)$ and are transversal to $\Lambda_1(L_0)$, 
then, by \eqref{eq:util}, $\mu_{L_0}$ is computed as
the difference between the number of positively oriented intersections
and the number of negatively oriented intersections of $l$ with
$\Lambda_1(L_0)$.

Using \eqref{eq:charPsi} and \eqref{eq:defellt}, a simple computation
yields:
\begin{equation}\label{eq:elllinha}
\ell'(t)=\omega(X(t)\,\cdot\,,\,\cdot\,)\vert_{\ell(t)}.
\end{equation}
\begin{teo}\label{thm:maslov=focal}
Let $(X,\ell_0)$ be a set of data for the symplectic differential
problem.
For $\varepsilon>0$ sufficiently small, there are no focal instants in
$]a,a+\varepsilon]$. If $t_0\in\,]a,b]$ is a nondegenerate focal instant, 
then $t_0$ is an  {\em isolated\/}  focal instant. If this nondegeneracy
condition is satisfied by all focal instants and if $t=b$
is not  focal, then the following equality holds:
\begin{equation}\label{eq:maslov=focal}
\maslov(X,\ell_0)=\mathrm i_{\mathrm{foc}}(X,\ell_0).
\end{equation}
\end{teo}
\begin{proof}
For all $t\in[a,b]$, the projection onto the second coordinate gives
an identification between $\ell(t)\cap L_0$ and the space $\mathbb V[t]^o$
(recall formula~\eqref{eq:VTo}). From \eqref{eq:elllinha}, it
follows easily that this identification carries the restriction of
$\ell'(t)$ to the restriction of $B(t)$. The conclusion follows then
from Proposition~\ref{thm:compmaslov}, observing that, by Assumption~\ref{thm:nondeg},
$B(a)$ is non degenerate on $P^o=\mathbb V[a]^o$. \end{proof}
\begin{rem}\label{thm:remuniform}
Suppose that we are given a family $\{X_\delta\}$  of coefficient matrices
for  symplectic differential systems 
such that $(t,\delta)\mapsto X_\delta(t)$ is continuous and $\delta$ runs
on a compact topological space. Then, it is not hard to show that we can find $\varepsilon>0$
independent on $\delta$ such that there are no $(X_\delta,\ell_0)$-focal instants
in $]a,a+\varepsilon]$. 
\end{rem}
\begin{rem}\label{thm:remstability}
By the homotopical invariance of $\mu_{L_0}$ and by 
Remark~\ref{thm:remuniform}, it follows that
the Maslov index of $(X,\ell_0)$ is stable by continuous deformations of 
the coefficient matrix $X$,
as long as the instant $t=b$ remains non focal during the deformation.
In particular, if $t=b$ is not focal, the Maslov index is stable by uniformly 
small perturbations of $X$.
\end{rem}
\subsection{On the nondegeneracy condition}
\label{sub:nondegeneracy}
We have seen in the previous subsection that the equality between the Maslov and
the focal index of a pair $(X,\ell_0)$ holds under
the assumption of nondegeneracy for the $(X,\ell_0)$-focal instants. 
We will call {\em nondegenerate\/} a pair
$(X,\ell_0)$ for which all the focal instants are nondegenerate. The aim of this section is
to prove that this condition is {\em generic}, i.e., that every pair $(X,\ell_0)$ can be
uniformly approximated by nondegenerate pairs.
\begin{prop}\label{thm:nondegeneracy}
Let $(X,\ell_0)$ be a set of data for the symplectic differential problem.
There exists a sequence $X_k:[a,b]\to\mathrm{sp}(2n,\R)$ of smooth curves
such that $X_k$ tends to $X$ uniformly as $k\to\infty$, and $(X_k,\ell_0)$ is a
nondegenerate set of data for the symplectic differential problem
for all $k$. Moreover, if $t=b$ is not $(X,\ell_0)$-focal, then, for
$k$ sufficiently large, $t=b$ is not $(X_k,\ell_0)$-focal.
\end{prop}
\begin{proof}
Since $X$ can be uniformly approximated by smooth curves, there
is no loss of generality in assuming that $X$ is smooth. 
The first step is to prove that $\ell$ can be approximated in the
$C^1$-topology by smooth curves of Lagrangians starting at $\ell_0$ having intersections
with $\Lambda_{\ge1}(L_0)$ only in $\Lambda_1(L_0)$ and transversally.

Let $\varepsilon>0$ be such that $\ell\vert_{]a,a+\varepsilon]}$ 
does not intersect $\Lambda_{\ge1}(L_0)$. It follows from \cite[2.1.\ Transversality
Theorem (a)]{Hirsch} that we can find a sequence $\ell_k:[a+\varepsilon,b]\to\Lambda$
of smooth curves that are transverse to $\Lambda_r(L_0)$ for all
$r=1,\ldots,n$ and such that $\ell_k$ converges to $\ell\vert_{[a+\varepsilon,b]}$
in the $C^1$-topology. Since $\Lambda_r(L_0)$ has codimension greater than one for
$r\ge2$, transversality to $\Lambda_r(L_0)$ in facts implies that $\ell_k$ does
not intercept $\Lambda_r(L_0)$. Hence, $\ell_k$ has only transverse intersections
with $\Lambda_1(L_0)$. Clearly, each $\ell_k$ can be extended to a smooth curve
on $[a,b]$ such that $\ell_k=\ell$ in $[a,a+\frac\varepsilon2]$ and such that $\ell_k$
converges to $\ell$ on $[a,b]$ in the $C^1$-topology. For $k$ sufficiently large,
it will also follow that $\ell_k$ does not intercept $\Lambda_{\ge1}(L_0)$ in $]a,a+\varepsilon]$.

It's now possible\footnote{%
This can be done, for instance, in the following way. Consider an arbitrary 
principal connection on the principal bundle $\mathrm{Sp}(2n,\R)\to\Lambda$ and let
$\hat\Psi_k$ be the parallel lifting of $\ell_k$ with $\hat\Psi_k(a)=\mathrm{Id}$;
then, $\Psi_k$ tends in the $C^1$-topology to the parallel lifting $\hat\Psi$ of $\ell$
starting at the identity. Define $\Psi_k=\hat\Psi_k\hat\Psi^{-1}\Psi$.} 
to construct a sequence $\Psi_k:[a,b]\to\mathrm{Sp}(2n,\R)$ of smooth
curves converging to $\Psi$ in the $C^1$-topology such that $\Psi_k(a)=\mathrm{Id}$ for all
$k$ and such that each $\Psi_k$ projects into $\ell_k$ by evaluation at $\ell_0$.

We now set, for each $k$, $X_k=\Psi_k^{-1}\Psi'_k$. Obviously $X_k$ is a smooth curve in
$\mathrm{sp}(2n,\R)$ and $X_k$ tends uniformly to $X$. Moreover, denoting by $B_k$ the
upper-right $n\times n$ block of $X_k$, we see that $B_k$ is invertible for $k$
sufficiently large, since $B_k$ tends uniformly to $B$ and $B$ is invertible. Similarly,
since $B(a)$ is nondegenerate on $P^o$ (recall assumption~\ref{thm:nondeg}), $B_k(a)$ is
nondegenerate on
$P^o$ for $k$ sufficiently large. Therefore $(X_k,\ell_0)$ is a set of data for
the symplectic differential problem, and obviously $\ell_k$ is its associated curve of
Lagrangians. It follows that $(X_k,\ell_0)$ is nondegenerate.

Finally, if $t=b$ is not $(X,\ell_0)$-focal then $\ell(b)\in\Lambda_0(L_0)$. Since
$\Lambda_0(L_0)$ is open and $\ell_k(b)$ tends to $\ell(b)$ we see that $t=b$ is not
$(X_k,\ell_0)$-focal for $k$ sufficiently large. 
\end{proof}
We conclude with the observation that, by Remark~\ref{thm:remstability},
in Proposition~\ref{thm:nondegeneracy} we have:
\[\maslov(X_k,\ell_0)=\maslov(X,\ell_0),\]
for $k$ sufficiently large.
\subsection{The index form}
\label{sub:indexform}
We denote by $L^2([a,b];\R^n)$ the Hilbert space of square integrable
$\R^n$-valued functions on $[a,b]$, and, for $j\ge1$, by $H^j([a,b];\R^n)$ the Sobolev
space of functions of class $C^{j-1}$, with $(j-1)$-th derivative
absolutely continuous  and with square integrable $j$-th derivative. 
We also denote by
$H^1_0([a,b];\R^n)$ the subspace of
$H^1([a,b];\R^n)$ consisting of those functions vanishing at $a$ and $b$; for a given
subspace $P\subset\R^n$ let $H^1_P([a,b];\R^n)$ be the subspace of $H^1([a,b];\R^n)$
consisting of those maps $v$ such that $v(a)\in P$ and $v(b)=0$.

Let $(X,\ell_0)$ be a fixed set of data for the symplectic differential problem
in $\R^n$; recall that there is a subspace $P\subset\R^n$ and a symmetric
bilinear form $S:P\times P\to\R$ canonically associated to $\ell_0$
as in \eqref{eq:bijection}. We will take into consideration the Hilbert space
$\mathcal H=H^1_P([a,b];\R^n)$; let $I$ be the following bounded symmetric bilinear form
on $\mathcal H$:
\begin{equation}\label{eq:defindexform}
\begin{split}
I(v,w)=&\int_a^b\Big[B(\alpha_v,\alpha_w)+C(v,w)\Big]\;\mathrm dt-S(v(a),w(a))=\\
=&\int_a^b\Big[B^{-1}(v'-Av,w'-Aw)+C(v,w)\Big]\;\mathrm dt-S(v(a),w(a)).
\end{split}
\end{equation}
We will call $I$ the {\em index form\/} associated to the pair $(X,\ell_0)$.
For instance, if the symplectic system comes from the Morse--Sturm equation
\eqref{eq:MS}, then $I$ becomes:
\begin{equation}\label{eq:indexformMS}
I(v,w)=\int_a^b\Big[g(v',w')+g(Rv,w)\Big]\;\mathrm dt-S(v(a),w(a));
\end{equation}
this is the classical index form for Morse--Sturm systems.

Integration by parts in \eqref{eq:defindexform} and the Fundamental Lemma of Calculus of
Variations show that:
\begin{equation}\label{eq:kerI}
\mathrm{Ker}(I)=\big\{v\in\mathbb V:v(b)=0\big\},
\end{equation}
where $\mathbb V$ is defined in \eqref{eq:spaces}.
\begin{rem}\label{thm:remposdef}
If $S$ is negative semidefinite, $B$ is positive definite and $C$ is positive semidefinite,
then  $I$ is positive definite in $\mathcal H$.
By considering restrictions to the interval $[a,t]$, it follows
easily from \eqref{eq:kerI} that the pair $(X,\ell_0)$ does not have focal instants in $]\,a,b]$.
Analogous results hold if $S$ is positive semidefinite, $B$ is negative definite
and $C$ is negative semidefinite.
\end{rem}
We will see in Subsection~\ref{sub:indextheorem} that if $B$ is positive
definite, then the index $n_-(I)$ in $\mathcal H$ is finite, and this number is
related to the focal instants of the symplectic differential problem.
However, for a non positive $B$, the index of $I$ in $\mathcal H$ is infinite:
\begin{prop}\label{thm:indinf}
The index of $I$ in $\mathcal H$ is infinite, unless $B$ is positive
definite.
\end{prop}
\begin{proof}
First, if $v$ is an $X$-solution
and $f\in H^1_0([a,b];\R)$, then, using \eqref{eq:alphafv},
we compute:
\begin{equation}\label{eq:interes}
I(fv,fv)=\int_a^b (f')^2B^{-1}(v,v)\;\mathrm dt.
\end{equation}
Let $t_0\in\,]a,b[$ and let $v$ be any $X$-solution 
such that $B(t_0)^{-1}(v(t_0),v(t_0))<0$, so that $B^{-1}(v,v)<0$
on $[t_0-\varepsilon,t_0+\varepsilon]$ for some $\varepsilon>0$. 
It follows from \eqref{eq:interes} that $I$ is negative definite on the
infinite dimensional space of maps $fv$, where $f$ is in $H_0^1([a,b];\R)$
and has support in $[t_0-\varepsilon,t_0+\varepsilon]$.
\end{proof}
One of the main goals of the paper is to determine a subspace of
$\mathcal H$ on which the index of $I$ is finite and that carries
the relevant information about the focal instants of the symplectic
problem.
\subsection{A subspace where the index is finite}
\label{sub:K}
We want to consider the case that the bilinear form $B$ is not necessarily 
positive definite; the idea is to {\em
factor out\/} an infinite dimensional subspace of $\mathcal H$ on which $I$ is negative
definite, and such that the kernel of $I$ remains unchanged.

Let's assume that the bilinear form $B$ has index $k$ in ${\R^n}^*$:
\begin{equation}\label{eq:indexB}
n_-(B)=n_-(B^{-1})=k.
\end{equation}
Let $\mathcal D=\{\mathcal D_t\}_{t\in[a,b]}$ be a  distribution
of subspaces $\mathcal D_t\subset\R^n$, with $\mathrm{dim}(\mathcal D_t)=k$;
assume that $B(t)^{-1}$ is negative definite on $\mathcal D_t$ for each $t$, and
that $\mathcal D_t$ has a $C^2$-dependence on $t$, in the sense that there
exists a  basis $Y_1(t),\ldots,Y_k(t)$ of $\mathcal D_t$, where each $Y_i$ is
of class $C^2$. We introduce the following
closed subspace of $\mathcal H$:
\begin{equation}\label{eq:newK}
\begin{split}
\mathcal K=\big\{v\in\mathcal H:&\;\alpha_v(Y_i)\in H^1([a,b];\R)\
\text{and}
\\&
\alpha_v(Y_i)' =B(\alpha_v,\alpha_{Y_i})+C(v,Y_i),\ \ \forall\,i=1,\ldots,k\big\}.
\end{split}
\end{equation}
It is not hard to see that $\mathcal K $ depends {\em only\/} on the distribution
of spaces $\mathcal D$, and not on the choice of the frame $Y_1,\ldots, Y_k$.
Roughly speaking, the space $\mathcal K$ is to be interpreted as the
space of those $v\in\mathcal H$ that are {\em $X$-solutions 
along $\mathcal D$}; namely, if $v\in\mathcal H$ is of class $C^2$, then
an easy computation shows that $v\in\mathcal K$ if and only if:
\[\alpha_v'(Y_i)=\big(Cv-A^*\alpha_v\big)(Y_i),\quad i=1,\ldots,k.\]
In particular, by \eqref{eq:kerI} we have:
\begin{equation}\label{eq:kerIinK}
\mathrm{Ker}(I)\subset\mathcal K.
\end{equation}
We associate to 
the choice of $Y_1,\ldots, Y_k$ the following  symplectic differential system in $\R^k$:
\begin{equation}\label{eq:symplassoc}
\left\{\begin{array}{l}f'=-\mathcal B^{-1}\mathcal C_{\mathrm a}f+\mathcal B^{-1}\varphi;\\
\varphi'=(\mathcal I-\mathcal C_{\mathrm s}'+\mathcal C_{\mathrm a}\mathcal B^{-1}
\mathcal C_{\mathrm a})f-\mathcal C_{\mathrm a}\mathcal B^{-1}\varphi,
\end{array}\right.\quad t\in [a,b],\ f(t)\in\R^k,\ \varphi(t)\in{\R^k}^*,
\end{equation}
where $\mathcal B$, $\mathcal I$
are bilinear forms in $\R^k$, and $\mathcal C$,
$\mathcal C_{\mathrm a}$, $\mathcal C_{\mathrm
s}$ are linear maps from $\R^k$ to ${\R^k}^*$, whose
matrices in the canonical basis are given by:
\begin{equation}
\label{eq:defcalobj}
\begin{split}
& \mathcal B_{ij}=B^{-1}(Y_i,Y_j),\quad \mathcal C_{ij}=\alpha_{Y_j}(Y_i),
\quad \mathcal C_{\mathrm s}=\frac12(\mathcal C+\mathcal C^*),\\
&\mathcal C_{\mathrm a}=\frac12(\mathcal C-\mathcal C^*),\quad
\mathcal I_{ij}=B(\alpha_{Y_i},\alpha_{Y_j})+C(Y_i,Y_j).
\end{split}
\end{equation}
We call \eqref{eq:symplassoc} the {\em reduced
symplectic system\/} associated to \eqref{eq:sympl} with respect
to $Y_1,\ldots,Y_k$.

We will make  the following Assumption  for Subsections~\ref{sub:K}, \ref{sub:indextheorem}
and
\ref{sub:variable}:
\begin{assum}\label{thm:assK}
We assume that the symplectic differential system \eqref{eq:symplassoc}
with initial condition $f(a)=0$ has no focal instants in $]a,b]$, i.e.,
if $f$ is a non zero solution of \eqref{eq:symplassoc} vanishing at $t=a$, then
$f$ does not vanish in $]a,b]$.
\end{assum}
Observe that the coefficient $\mathcal B^{-1}$ in  \eqref{eq:symplassoc}
is negative definite, and in Remark~\ref{thm:remposdef} we gave a criterion
for the absence of focal instants: if the symmetric  bilinear form
$\mathcal I-\mathcal C_{\mathrm s}'+\mathcal C_{\mathrm a}\mathcal B^{-1}
\mathcal C_{\mathrm a}$ is negative semidefinite, then Assumption~\ref{thm:assK}
is satisfied (see also Remark~\ref{thm:remfacilitar} ahead).
\begin{rem}\label{thm:remdoesnotdepend}
Although the definition of the system \eqref{eq:symplassoc}
depends on the choice of a basis $\{Y_1,\ldots,Y_k\}$ of $\mathcal D$, 
the reduced symplectic
systems associated to different choices of frames of $\mathcal D$ are {\em isomorphic}, in a sense
that will be clarified in Subsection~\ref{sub:equivalence} (see
Proposition~\ref{thm:reducedisomorphic}). In particular, the validity of Assumption~\ref{thm:assK}
does not depend on the choice of $Y_1,\ldots,Y_k$.
\end{rem}
\noindent The motivation for the definition of the reduced symplectic system
is given by Lemma~\ref{thm:KinterS} below.

The reader can keep in mind
the following examples where Assumption~\ref{thm:assK} holds.

\begin{example}\label{thm:example1}
If $k=0$, then there is no need of any assumption. If $k\ge1$ and there exist
pointwise linearly independent $X$-solutions $Y_1,\ldots,Y_k$ such that:
\begin{enumerate}
\item\label{itm:hp1} $B^{-1}$ is negative definite on the span
of $Y_1,\ldots,Y_k$ for all $t$;
\item\label{itm:hp2} the matrix $\alpha_{Y_i}(Y_j)$ is symmetric for all $t$,
\end{enumerate}
then the right hand side of the second equation of \eqref{eq:symplassoc}
vanishes identically, and therefore Assumption~\ref{thm:assK}
is satisfied. In this case, it is not hard to see that the space $\mathcal K$ defined in
\eqref{eq:newK} is   given by:
\[\mathcal K=\big\{v\in\mathcal H:\alpha_v(Y_i)-\alpha_{Y_i}(v)\ \text{is constant}
\ \ \forall\,i=1,\ldots,k\big\}.\]
For instance, if $k=1$ then, in order to satisfy Assumption~\ref{thm:assK}, 
one only needs to determine an $X$-solution
$Y$ such that $B^{-1}(Y,Y)$ is negative.
\end{example}
\begin{example}\label{thm:example2}
If the original symplectic system is a Morse--Sturm system of the
form \eqref{eq:MS}, then an example where
Assumption~\ref{thm:assK} holds is obtained when it is possible to find a constant
$k$-dimensional subspace $\mathcal D$ on which  $B(t)^{-1}$ 
is negative definite and $C(t)$ is negative semi-definite
for all $t$. In this case, we can take constant vector fields $Y_1,\ldots,Y_k$
as a basis for~$\mathcal D$.
\end{example}

We now discuss some properties of the space $\mathcal K$.
\begin{rem}\label{thm:remcompact}
Let $E$ be a closed subspace of $H^1([a,b];\R^n)$, suppose that $B\in\mathrm B(E,\R)$
is a bilinear form. Let $E_0$ denote the normed vector space obtained by considering
the $C^0$-topology on $E$ and assume that $B$ is continuous in $E_0\times E$.
It follows that the linear operator $T\in\mathcal L(E)$ which represents
$B$ is {\em compact}, since the inclusion of $H^1$ in $C^0$ is compact.
The same conclusion holds if we assume that $B$ is continuous in $E\times E_0$, since
the adjoint  of a compact operator in a Hilbert space is also compact.
\end{rem}
We are now ready for the following:
\begin{lem}\label{thm:P+K}
The restriction of $I$ to $\mathcal K$ is represented by a compact perturbation
of a positive isomorphism.
\end{lem}
\begin{proof}
Consider the bilinear form $I_0(v,w)=\int_a^bB^{-1}(v',w')\;\mathrm dt$;
recalling \eqref{eq:defindexform}, from Remark~\ref{thm:remcompact} it 
follows easily that $I-I_0$ is represented by a compact operator on $\mathcal K$.
It remains to prove that $I_0$ is represented by a compact perturbation of a positive
isomorphism of $\mathcal K$. For each $t$, we define a positive definite inner product
$r_t$ on $\R^n$ by setting 
\[r_t(x,y)=B(t)^{-1}(x,y)-2\,B(t)^{-1}(\pi(t)x,\pi(t)y),\]
where $\pi(t):\R^n\to\mathcal D(t)$ is the orthogonal projection with respect
to $B(t)^{-1}$. We have:
\begin{equation}\label{eq:r-B}
B^{-1}(v,w)=r(v,w)-2\,r(\pi(v),\pi(w)),
\end{equation}
and it follows:
\[I_0(v,w)=\int_a^b\Big[r(v',w')-2\,r(\pi(v'),\pi(w'))\Big]\;\mathrm dt.\]
Clearly, the integral of the first term above gives a Hilbert space inner
product in $\mathcal K$, and it is therefore represented by a positive
isomorphism of $\mathcal K$. Using again Remark~\ref{thm:remcompact}, to conclude
the proof it suffices to show that the linear map $\mathcal K\ni v\mapsto\pi(v')
\in L^2([a,b];\R^n)$ is continuous in the $C^0$-topology of $\mathcal K$.
To this aim, it clearly suffices to show that the maps $v\mapsto \alpha_v(Y_i)$
are $C^0$-continuous for all $i$. This follows from the fact that, for $v\in\mathcal K$
the quantity $c_i(v)$ defined by:
\begin{equation}\label{eq:expand}
c_i(v)=\alpha_v(Y_i)-\int_a^t\Big[B(\alpha_v,\alpha_{Y_i})+C(v,Y_i)\Big]\;\mathrm
ds 
\end{equation}
is constant (see formula \eqref{eq:newK}). Using \eqref{eq:defalphav}
and integration by parts, we see that the integral in \eqref{eq:expand}
is continuous in $v$ with respect to the $C^0$-topology. Integrating
\eqref{eq:expand} on $[a,b]$, we see that the functional $\mathcal K\ni v\mapsto c_i(v)$
is $C^0$-continuous, which concludes the proof.
\end{proof}
\begin{cor}\label{thm:finiteindex}
$I$ has finite index on $\mathcal K$.\qed
\end{cor}
We now describe $\mathcal K$ as the kernel of a bounded linear
map $F:\mathcal H\to L^2([a,b];{\R^k}^*)/\mathfrak C$, where $\mathfrak C$
is the subspace of constant maps. Let $F=(F^1,\ldots,F^k)$ be defined by:
\begin{equation}\label{eq:defFi}
F^i(v)(t)=\alpha_v(t)(Y_i(t))-\int_a^t\Big[B(\alpha_v,\alpha_{Y_i})+C(v,Y_i)\Big]\;\mathrm
ds +\ \text{constant}.
\end{equation}
From \eqref{eq:newK} it follows that $\mathcal K=\mathrm{Ker}(F)$;
moreover, $F$ is clearly continuous in $\mathcal H$, and therefore $\mathcal K$
is closed. Next, we introduce a closed subspace $\mathcal S$, that will be shown
to be a closed complement of $\mathcal K$ in $\mathcal H$:
\begin{equation}\label{eq:defS}
\mathcal S=\big\{v\in H^1_0([a,b];\R^n):\ v(t)\in\mathcal D(t),\ \forall\,t\in[a,b]\big\}.
\end{equation}
\begin{lem}\label{thm:KinterS}
$\mathcal K\cap\mathcal S=\{0\}$.
\end{lem}
\begin{proof}
If $v=\sum_i f_iY_i\in\mathcal S$,  then $v\in\mathcal K$ iff
$f=(f_1,\ldots,f_k)$ is a solution of the reduced symplectic system
\eqref{eq:symplassoc}.
\end{proof}
\begin{lem}\label{thm:Fsobre}
The restriction of $F$ to $\mathcal S$ is onto
$L^2([a,b];{\R^k}^*)/\mathfrak C$.
\end{lem}
\begin{proof}
We identify $\mathcal S$ with $H^1_0([a,b];\R^k)$ by $v=\sum_if_iY_i\mapsto f=(f_1,\ldots,
f_k)$. For $v\in\mathcal S$, we can rewrite $F(v)$ with the help of \eqref{eq:defcalobj}
as:
\begin{equation}\label{eq:iso+comp}
F(v)(t)=\mathcal B(t)f'(t)+\mathcal C(t)f(t)-\int_a^t\Big[\mathcal
C^*f'+\mathcal If\Big]\;\mathrm ds +\text{constant}.
\end{equation}
Define $F_0(f)=\mathcal Bf'$; then, $F_0$ is a linear isomorphism between  the spaces
$H^1_0([a,b];\R^k)$ and $L^2([a,b];{\R^k}^*)/\mathfrak C$.
Using integration by parts, the formula \eqref{eq:iso+comp} shows that 
the linear operator $F-F_0$ is continuous in the $C^0$-topology, and it
is therefore compact. Then, $F$ is a Fredholm operator of index zero on 
$\mathcal S$, and the injectivity of $F$ on $\mathcal S$
(Lemma~\ref{thm:KinterS}) implies that $F$ is onto. 
\end{proof}
\begin{cor}\label{thm:H=K+S}
$\mathcal H=\mathcal K\oplus\mathcal S$.
\end{cor}
\begin{proof} 
It follows from Lemma~\ref{thm:KinterS}, Lemma~\ref{thm:Fsobre} and the fact
that $\mathcal K=\mathrm{Ker}(F)$.
\end{proof}
\begin{lem}\label{thm:SKorth}
$\mathcal K$ and $\mathcal S$ are $I$-orthogonal, i.e., $I(v,w)=0$ for all
$v\in\mathcal S$ and $w\in\mathcal K$.
\end{lem}
\begin{proof}
It is an easy computation.
\end{proof}
\begin{cor}\label{thm:kerIK}
$\mathrm{Ker}(I\vert_{\mathcal K})=\mathrm{Ker(I)}$.
\end{cor}
\begin{proof}
It follows easily from \eqref{eq:kerIinK}, Corollary~\ref{thm:H=K+S}
and Lemma~\ref{thm:SKorth}.
\end{proof}
\subsection{The index theorem}
\label{sub:indextheorem}
We have proven the finiteness of the index of $I$ in $\mathcal K$;
in this section we determine the value of this index in terms of
the Maslov index of the pair $(X,\ell_0)$. To this aim, we study the
{\em evolution\/} of the index function $i(t)$, $t\in]a,b]$, defined as
the index of the form $I_t$ on the space $\mathcal K_t\subset\mathcal H_t$, 
where $\mathcal H_t=H^1_P([a,t];\R^n)$ and $I_t$ and $\mathcal K_t$ are defined 
as in formulas \eqref{eq:defindexform} and \eqref{eq:newK} by replacing
$b$ with $t$ and $\mathcal H$ with $\mathcal H_t$. Similarly, one can define 
in an obvious way the {\em
objects\/} $F_t$ and $\mathcal S_t$ as in \eqref{eq:defFi} and \eqref{eq:defS};
clearly, all the results proven in Subsection~\ref{sub:K} remain valid
when the symplectic differential system is restricted to the
interval $[a,t]$. Observe that for this reason in Assumption~\ref{thm:assK}
we have required that the reduced symplectic system \eqref{eq:symplassoc} had no
focal instants in the entire interval $]a,b]$, and not just at the final
instant $t=b$.

We will use the isomorphisms $\Phi_t:\mathcal H\to\mathcal H_t$ defined by:
$\Phi_t (\hat v)=v$, where 
\begin{equation}\label{eq:Phit}
v(s)=\hat v(u_s),\quad u_s=a+\frac{b-a}{t-a}(s-a),\quad\forall\,s\in[a,t];
\end{equation}
and we get a family $\{\hat{\mathcal K}_t\}_{t\in]a,b]}$ of closed subspaces
of $\mathcal H$, a curve $\hat I:\,]a,b]\to\Bsym(\mathcal H,\R)$ of
symmetric bilinear forms and a curve $\hat F:\,]a,b]\to \mathcal L(\mathcal H,
L^2([a,b];{\R^k}^*)/\mathfrak C)$ of maps, defined by:
\[\hat{\mathcal K}_t=\Phi_t^{-1}(\mathcal K_t),\quad \hat
I_t=I(\Phi_t\,\cdot\,,\Phi_t\,\cdot\,),\quad\hat F_t=\Phi_t^{-1}\circ F_t\circ
\Phi_t.\]
We are also denoting by $\Phi_t$ the isomorphism from $L^2([a,b];{\R^k}^*)/\mathfrak C$
to $L^2([a,t];{\R^k}^*)/\mathfrak C$ defined by formula \eqref{eq:Phit}.
An explicit formula for $\hat I_t$ is given by:
\begin{equation}\label{eq:explhatIt}
\begin{split}
\hat I_t(\hat v,\hat w)=& \int_a^t\! B(s)^{-1}\left(\frac{b-a}{t-a}\hat v'(u_s)-
A(s)(\hat v(u_s)),\frac{b-a}{t-a}\hat w'(u_s)-
A(s)(\hat w(u_s))\right) \mathrm ds\\+& 
\int_a^t C(s)(\hat v(u_s),\hat w(u_s))\;\mathrm ds-S(\hat v(a),\hat w(a)).
\end{split}
\end{equation}
For $t\in\,]a,b]$, we set $\mathcal J_t=(t-a)\hat I_t$:
\begin{equation}\label{eq:Jcalt}
\begin{split}
\mathcal J_t(\hat v,\hat w)=&\, \int_a^bB(s_u)^{-1}((b-a)\hat v'(u)
-(t-a)A(s_u)\hat v(u),\hat w'(u))\;\mathrm du  \\
-\frac1{b-a}&\int_a^b B(s_u)^{-1}((b-a)\hat v'(u)
-(t-a)A(s_u)\hat v(u),(t-a)A(s_u)\hat w(u))\;\mathrm du\\
&+\frac1{b-a}\int_a^b(t-a)^2C(s_u)(\hat v(u),\hat w(u))\;\mathrm du-(t-a)S(\hat v(a),
\hat w(a)),
\end{split}
\end{equation}
where $s_u=a+\frac{t-a}{b-a}(u-a)$.
Setting $t=a$, formula \eqref{eq:Jcalt} defines $\mathcal J_a$ as:
\begin{equation}\label{eq:defJcala}
\mathcal J_a(\hat v,\hat w)=(b-a)\int_a^bB(a)^{-1}(\hat v'(u),\hat w'(u))\;\mathrm du.
\end{equation}
Obviously, for $t\in\,]a,b]$, we have
\begin{equation}\label{eq:defit}
i(t)=n_-(I_t\vert_{\mathcal K_t})=n_-(\hat
I_t\vert_{\hat{\mathcal K}_t}) =n_-(\mathcal J_t\vert_{\hat{\mathcal K}_t}).
\end{equation}
In order to study the evolution of the function $i$, we will make full
use of some results presented in reference~\cite{GMPT} concerning the
jumps of the index of a $C^1$-curve of bounded symmetric bilinear forms
restricted to a $C^1$-family of closed subspaces of a fixed Hilbert space.
For the reader's convenience, we recall some definitions and facts presented
in \cite{GMPT}.

\begin{defin}\label{thm:defC1subspaces}
Let $H$ be a Hilbert space, $I\subset\R$ an interval and $\{D_t\}_{t\in I}$
be a family of closed subspaces of $H$. We say that $\{D_t\}_{t\in I}$
is a {\em $C^1$-family\/} of subspaces if for all $t_0\in I$ there exists a $C^1$-curve
$\beta:\,]t_0-\varepsilon,t_0+\varepsilon\,[\,\cap\, I\to \mathcal L(H)$ 
and a closed subspace $\overline{D}\subset H$
such that $\beta(t)$ is an isomorphism  and $\beta(t)(D_t)=\overline{D}$
for all $t$.
\end{defin}
We have a method for producing $C^1$-families of closed subspaces:
\begin{prop}\label{thm:produce}
Let $I\subset\R$ be an interval, $H,\tilde{H}$ be Hilbert
spaces and $\mathcal F:I\to\mathcal L(H,\tilde{H})$ be a $C^1$-map
such that each $\mathcal F(t)$ is surjective. Then, the family
$D_t=\mathrm{Ker}(\mathcal F(t))$ is a $C^1$-family of closed subspaces of $H$.
\end{prop}
\begin{proof} See \cite[Lemma~2.9]{GMPT}.
\end{proof}
Here comes a method for computing the jumps of 
the index of a smooth curve $M(t)$ of symmetric bilinear forms 
represented by a compact perturbation of positive isomorphisms.
\begin{prop}\label{thm:HSelementary}
Let $ H$ be a real Hilbert space with inner product
$\langle\cdot,\cdot\rangle$, and let $M:[t_0,t_0+r]\to \mathrm
B_{\mathrm{sym}}( H,\R)$, $r>0$, be a map of class $C^1$. 
Let $\{D_t\}_{t\in[t_0,t_0+r]}$ be a $C^1$-family of closed subspaces
of $H$, and denote by $\overline M(t)$ the restriction of $M(t)$
to $D_t\times D_t$.
Assume that the following three hypotheses are satisfied:
\begin{enumerate}
\item\label{itm:hp2.5.1} $\overline M(t_0)$ is represented by a compact perturbation
of a positive isomorphism of $D_{t_0}$;
\item\label{itm:hp2.5.2} the restriction $\widetilde M$ of the derivative $M'(t_0)$ to 
$\mathrm{Ker}(\overline M(t_0))\times \mathrm{Ker}(\overline M(t_0))$ 
is non degenerate;
\item\label{itm:hp2.5.3} $\mathrm{Ker}(\overline M(t_0))\subseteq \mathrm{Ker}(M(t_0))$.
\end{enumerate}
Then, for
$t>t_0$ sufficiently close to $t_0$, $\overline M(t)$ is non degenerate, and we have:
\begin{equation}\label{eq:changen-}
n_-(\overline M(t))=n_-(\overline M(t_0))+n_-(\widetilde M),
\end{equation}
all the terms of the above equality being finite natural numbers.
\end{prop}
\begin{proof} See \cite[Proposition~2.5]{GMPT}.\end{proof}
We have the following two immediate corollaries of Proposition~\ref{thm:HSelementary}.
\begin{cor}\label{thm:jump}
Under the hypotheses of Proposition~\ref{thm:HSelementary}, if
$M$ and $D_t$ are defined and of class $C^1$ in a neighborhood
of $t_0$, then, for $\varepsilon>0$ sufficiently small, we have:
\begin{equation}\label{eq:jump}
n_-(\overline M(t_0+\varepsilon))-n_-(\overline M(t_0-\varepsilon))=
-\sgn(\widetilde M).
\end{equation}
\end{cor}
\begin{proof}
Apply twice Proposition~\ref{thm:HSelementary}, once to $M$ and again to a 
{\em backwards\/} reparameterization of $M$.
\end{proof}
\begin{cor}\label{thm:cornondeg}
Under the hypotheses of Proposition~\ref{thm:HSelementary}, if $\overline M(t_0)$
is nondegenerate, then $n_-(\overline M(t))$ is constant for $t$ near $t_0$.
\hfill\qed
\end{cor}

We now show how our setup fits into the assumptions of Proposition~\ref{thm:HSelementary}.

For $t\in\,]a,b]$, we set $\mathcal F_t=(t-a)\hat F_t$; explicitly, 
for $i=1,\ldots,k$, $\hat v\in\mathcal H$ and
$u\in[a,b]$, we have:
\begin{equation}\label{eq:explicitly}
\begin{split}
\mathcal F_t^i(\hat v)(u)=&\; B(s_u)^{-1}((b-a)\hat v(u)-(t-a)A(s_u)\hat v(u),Y_i(s_u))\\ &
-\frac{s_u-a}{u-a}\int_a^u \alpha_{Y_i}(r_x)((b-a)\hat v'(x)-(t-a)A(r_x)\hat v(x))
\;\mathrm dx\\&-\frac{s_u-a}{u-a}\int_a^u(t-a)
C(r_x)(\hat v(x),Y_i(r_x)) \;\mathrm dx+\text{constant},
\end{split}
\end{equation}
where $r_x=a+\frac{s_u-a}{u-a}(x-a)$.
Setting $t=a$, formula \eqref{eq:explicitly} defines $\mathcal F_a$ as:
\begin{equation}\label{eq:defFcala}
\mathcal F_a^i(\hat v)(u)=(b-a)B(a)^{-1}(\hat v'(u),Y_i(a))+\text{constant}.
\end{equation}
Obviously, $\mathrm{Ker}(\mathcal F_t)=\hat{\mathcal K}_t$ for $t\in\,]a,b]$;
we set $\hat{\mathcal K}_a=\mathrm{Ker}(\mathcal F_a)$, namely:
\begin{equation}\label{eq:hatKa}
\hat{\mathcal K}_a=\big\{\hat v\in\mathcal H:B(a)^{-1}(\hat
v'(u),Y_i(a))=\text{constant},\  i=1,\ldots,k\big\}.
\end{equation}
\begin{prop}\label{thm:tudoC1}
Assume that $C$ is a map of class $C^1$. Then, $\hat I$ is a map of class $C^1$
in $]a,b]$, $\mathcal J$ is of class $C^1$ in $[a,b]$ and $\{\hat{\mathcal K}_t\}_{t
\in[a,b]}$ is a $C^1$-family of closed subspaces of $\mathcal H$.
\end{prop}
\begin{proof}
By standard regularity arguments (see \cite[Lemma~2.3, Proposition~3.3 and
Lem\-ma~4.3]{GMPT}), formula \eqref{eq:Jcalt} shows that $\mathcal J$ is $C^1$
in $[a,b]$, which obviously implies that $\hat I$ is $C^1$ in $]a,b]$.
Similarly, formula \eqref{eq:explicitly} shows that $\mathcal F$ is of
class $C^1$ on $[a,b]$; from Lemma~\ref{thm:Fsobre} we deduce that
$\mathcal F_t$ is surjective for $t\in\,]a,b]$. The surjectivity of
$\mathcal F_a$ is immediately established from its definition \eqref{eq:defFcala}.
The regularity of the family $\{\hat{\mathcal K}_t\}_{t
\in[a,b]}$ follows then from Proposition~\ref{thm:produce}.
\end{proof}
The last ingredient missing for applying Proposition~\ref{thm:HSelementary}
to the curve $\hat I$ is the computation of the derivative of $\hat I$ on its
kernel. This is done in the following:
\begin{lem}\label{eq:compderivative}
For each $t\in\,]a,b]$, the map $\sigma_t:\mathrm{Ker}(I_t)\to
\mathbb V[t]^o$ given by $\sigma_t(v)=\alpha_v(t)$ is an isomorphism.
Moreover, if $C$ is of class $C^1$, the isomorphism $\sigma_t\circ\Phi_t$ 
carries the restriction of the derivative $\hat I'_t$ to $\mathrm{Ker}(\hat I_t)$
into the restriction of $-B(t)$ to $\mathbb V[t]^o$.
\end{lem}
\begin{proof}
The fact that $\sigma_t$ is an isomorphism follows easily from 
formulas \eqref{eq:VTo} and \eqref{eq:kerI}.  

Let $\hat v,\hat w\in\mathrm{Ker}(\hat I_t)$; set $v=\Phi_t(\hat v)$ and $w=\Phi_t(\hat w)$,
so that $v,w\in\mathrm{Ker}(I_t)$. From \eqref{eq:explhatIt}, considering that
$v(t)=w(t)=0$ we compute:
\begin{equation}\label{eq:chata}
\begin{split}
\frac{\mathrm d}{\mathrm dt}&\,\hat I_t(\hat v,\hat w)= \\
&B(t)^{-1}(v'(t),w'(t))-\frac1{t-a}\int_a^t\Big[2B^{-1}(v',w')-B^{-1}(v',Aw)
-B^{-1}(Av,w')\Big]\;\mathrm ds\\ -&
 \frac1{t-a}\int_a^t\!\!(s-a)\Big[B^{-1}(v'',w')+B^{-1}(v',w'')-B^{-1}(v'',Aw)-B^{-1}(v',Aw')
\Big]\mathrm ds\\ &-\frac1{t-a}\int_a^t(s-a)\Big[
-B^{-1}(Av',w')-B^{-1}(Av,w'')+B^{-1}(Av',Aw)\Big]\;\mathrm ds\\
&-\frac1{t-a}\int_a^t(s-a)\Big[ B^{-1}(Av,Aw') +C(v',w)+C(v,w')\Big]\;\mathrm ds,
\end{split}
\end{equation}
where all the functions inside the integrals above are meant to be evaluated at $s$.
We use the equation  \eqref{eq:secord} satisfied by $v$ and $w$ to eliminate
the two terms containing the bilinear form $C$ in \eqref{eq:chata}, and we obtain:
\begin{equation}\label{eq:razoavel}
\begin{split}
\frac{\mathrm d}{\mathrm dt}\hat I_t(\hat v,\hat w)=&\;B(t)^{-1}(v'(t),w'(t))\\&-
\frac1{t-a}\int_a^t\frac{\mathrm d}{\mathrm ds}\,\Big[ (s-a) \big(
B^{-1}(v',w'-Aw)+B^{-1}(v'-Av,w')\big)\Big]\;\mathrm ds=\\ &
=-B(t)^{-1}(v'(t),w'(t))=-B(t)(\alpha_v(t),\alpha_w(t)).
\end{split}
\end{equation}
This concludes the proof.
\end{proof}
We now consider the case $t=a$.
\begin{lem}\label{thm:indJa}
The restriction of the symmetric bilinear form $\mathcal J_a$ to $\hat{\mathcal K}_a$
is represented by a compact perturbation of a positive isomorphism. Moreover, it
is non degenerate, and its index equals the index of the restriction of $B(a)^{-1}$
to $P$.
\end{lem}
\begin{proof}
The fact that the restriction of $\mathcal J_a$ to $\hat{\mathcal K}_a$
is represented by a compact perturbation of a positive isomorphism
follows by arguments similar to those used in the proof of Lemma~\ref{thm:P+K}.

Write $P=P_+\oplus P_-$, where $B(a)^{-1}$ is positive definite on
$P_+$ and negative definite on $P_-$; this is possible because, by
Assumption~\ref{thm:nondeg},
$B(a)^{-1}$ is nondegenerate on $P$. We now define the following
subspaces of $\hat{\mathcal K}_a$:
\begin{equation}\label{eq:K+K-}
\begin{split}
&\hat{\mathcal K}_+=\big\{\hat v\in\hat{\mathcal K}_a:\hat v(a)\in P_+\big\},\\
&\hat{\mathcal K}_-=\big\{\hat v:\hat v\ \text{is an affine
function}, \ \hat v(a)\in P_-,\ \hat v(b)=0\big\}.
\end{split}
\end{equation}
It is easily seen that $\hat{\mathcal K}_a=\hat{\mathcal K}_+\oplus\hat{\mathcal K}_-$;
we claim that $\mathcal J_a$ is positive definite on $\hat{\mathcal K}_+$
and negative definite on $\hat{\mathcal K}_-$. Namely, for $v_0\in P_-$,
$v_0\ne0$, and $\hat v\in \hat{\mathcal K}_-$ of the form
$\hat v(u)=v_0(u-a)$, from \eqref{eq:defJcala} it is easily computed 
$\mathcal J_a(\hat v,\hat v)=(b-a)^2B(a)^{-1}(v_0,v_0)<0$. Suppose now
$\hat v\in \hat{\mathcal K}_+$. Using the Cauchy--Schwarz inequality componentwise, 
it is easily proven the following inequality for any integrable function
$z:[a,b]\to\R^n$:
\begin{equation}\label{eq:Jensen}
r_a\left(\int_a^bz,\int_a^bz\right)\le (b-a)\int_a^br_a(z,z),
\end{equation}
where the equality holds if and only if $z$ is constant almost everywhere.
Applying \eqref{eq:Jensen} to $\hat v'$, we get:
\begin{equation}
\begin{split}
\mathcal J_a(\hat v,\hat v)\ge &\;r_a\left(\int_a^b\hat v',\int_a^b\hat v'\right)-2r_a(\pi(a)(\hat
v(a)),
\pi(a)(\hat v(a)))=\\&=B(a)^{-1}(\hat v(a),\hat v(a))\ge 0,
\end{split}
\end{equation}
where equality between the first and the last term above occurs if and only if
$\hat v$ is affine and $\hat v(a)=0$, i.e., if and only if $\hat v=0$.
This proves the claim; it follows easily that $\mathcal J_a$ is nondegenerate
on $\hat{\mathcal K}_a$ and that its index is equal to $\mathrm{dim}(\hat{\mathcal K}_-)
=\mathrm{dim}(P_-)=n_-(B(a)^{-1}\vert_P)$.
\end{proof}
We are ready to prove the index theorem for nondegenerate pairs $(X,\ell_0)$
(recall Subsection~\ref{sub:nondegeneracy}).
\begin{lem}\label{thm:Indexnondeg}
Let $(X,\ell_0)$ be a nondegenerate pair, with $X$ of class $C^1$.
Assume that $t=b$ is not a focal instant. Then, 
\begin{equation}\label{eq:indexth}
n_-(I\vert_{\mathcal K})=n_-(B(a)^{-1}\vert_P)+\mathrm i_{\mathrm{foc}}(X,\ell_0).
\end{equation}
\end{lem}
\begin{proof}
By Theorem~\ref{thm:maslov=focal}, the number of focal instants is finite.
Using \eqref{eq:kerI} and Corollary~\ref{thm:kerIK}, we see that $\hat I_t$ is
nondegenerate on $\hat{\mathcal K}_t$ when $t$ is not focal.
Applying Corollary~\ref{thm:cornondeg}, where Lemma~\ref{thm:P+K} is being used, 
the function $i(t)$ defined in \eqref{eq:defit} is piecewise constant, more precisely,
it is constant on every interval that does not contain focal instants.
Using Corollary~\ref{thm:jump} and Lemma~\ref{eq:compderivative}, the jump
of $i(t)$ at a focal instant $t\in\,]a,b[$ is equal to $\sgn(t)$. Finally,
by Corollary~\ref{thm:cornondeg} and Lemma~\ref{thm:indJa}, $i(t)=n_-(B(a)^{-1}\vert_P)$
for $t$ sufficiently close to $a$. This concludes the proof.
\end{proof}
We can now prove the aimed index theorem, which we state in a complete
form for future reference:
\begin{teo}[Index Theorem]\label{thm:indextheorem}
Let $(X,\ell_0)$ be a set of data for the symplectic differential problem
in $\R^n$, with $A,B$ of class $C^1$, $C$ continuous, $B(t)$ invertible
for all $t$, and let $\{\mathcal D_t\}_{t\in[a,b]}$ be a $C^2$-family of
$k$-dimensional subspaces of $\R^n$, where $k=n_-(B(t))$, and
$B(t)^{-1}$ is negative definite on $\mathcal D_t$ for all $t$.
Suppose that Assumption~\ref{thm:nondeg} and Assumption~\ref{thm:assK} are satisfied,
and that $t=b$ is not a focal instant. Define  $\mathcal K$ as in \eqref{eq:newK};
then, the index of $I$ on $\mathcal K$ is finite, and the following equality holds:
\begin{equation}\label{eq:equalityidex}
n_-(I\vert_{\mathcal K})=n_-(B(a)^{-1}\vert_P)+\maslov(X,\ell_0).
\end{equation}
\end{teo}
\begin{proof}
Using Theorem~\ref{thm:maslov=focal} and Lemma~\ref{thm:Indexnondeg}, the theorem
holds if $C$ is a map of class $C^1$ and if $(X,\ell_0)$ is nondegenerate.
The conclusion follows from Proposition~\ref{thm:nondegeneracy}, observing
that all the numbers in formula \eqref{eq:equalityidex} are stable
by uniformly small perturbations of the data.
\end{proof}
\begin{rem}\label{thm:remfinaldeg}
If $C$ is of class $C^1$, the result of Theorem~\ref{thm:indextheorem} can be extended to the case that
$t=b$ is a nondegenerate focal instant.
In this case, $t=b$ is an isolated focal instant, and one can define the Maslov
index of the pair $(X,\ell_0)$, in analogy with formula~\eqref{eq:defmaslov}, as:
\begin{equation}\label{eq:defmaslovbfocal}
\maslov=\mu_{L_0}(\ell\vert_{[a+\varepsilon,b-\varepsilon]}),
\end{equation}
where $\varepsilon>0$ is small enough. In this case, another application of
Proposition~\ref{thm:HSelementary} around $t=b$ gives us the following equality:
\begin{equation}\label{eq:equalityidexbfocal}
n_-(I\vert_{\mathcal
K})=n_-(B(a)^{-1}\vert_P)-n_-(B(b)^{-1}\vert_{\mathbb V[b]^o})+\maslov(X,\ell_0).
\end{equation}
\end{rem}
\subsection{The case of a variable endpoint}
\label{sub:variable}
Motivated by the geometric problem of solutions of a Hamiltonian
system with both endpoints variable, we will now consider the following
setup.

Given a pair $(X,\ell_0)$ of data for the symplectic differential
problem \eqref{eq:sympl}, we will additionally consider a Lagrangian
subspace $\ell_1$ of $(R^n\oplus{\R^n}^*,\omega)$, to which we associate
a pair $(Q,S_Q)$ consisting of a subspace $Q$ of $\R^n$ and
a symmetric bilinear form $S_Q$ on $Q$ as in \eqref{eq:bijection}.
Given the triple
$(X,\ell_0,\ell_1)$, we will consider a Hilbert space $\mathcal H^\#$ and
a bounded symmetric bilinear form $I^\#$ on $\mathcal H^\#$ as follows:
\begin{equation}\label{eq:defHsust}
\mathcal H^\#=\big\{v\in H^1([a,b];\R^n):v(a)\in P,\ v(b)\in Q\big\};
\end{equation}
and
\begin{equation}\label{eq:defIsust}
I^\#(v,w)=\int_a^b\Big[B(\alpha_v,\alpha_w)+C(v,w)\Big]\;\mathrm dt
+S_Q(v(b),w(b))-S(v(a),w(a)),
\end{equation}
with $\alpha_v,\alpha_w$ given by \eqref{eq:defalphav}.
Comparing with \eqref{eq:defindexform}, clearly $I^\#$ coincides with
$I$ on $\mathcal H$.

We want to extend the index theorem to this situation, and to this aim
we define a space $\mathcal K^\#$ in analogy with \eqref{eq:newK}:
\begin{equation}\label{eq:defKsust}
\begin{split}
\mathcal K^\#=\big\{v\in\mathcal H^\#:&\;\alpha_v(Y_i)\in H^1([a,b];\R)\
\text{and}
\\&
\alpha_v(Y_i)' =B(\alpha_v,\alpha_{Y_i})+C(v,Y_i),\ \ \forall\,i=1,\ldots,k\big\}.
\end{split}
\end{equation}
where $Y_1,\ldots,Y_k$ is a frame for $\mathcal D$..

Recalling the curve $\ell(t)$ of Lagrangians \eqref{eq:defellt} associated
to the pair $(X,\ell_0)$, we consider the pair $(\mathbb V_b,S_b)$ consisting of
a subspace $\mathbb V_b\subset\R^n$ and a symmetric bilinear form
$S_b$ on $\mathbb V_b$ associated to the Lagrangian $\ell(b)$; explicitly,
we have:
\begin{equation}\label{eq:Vb}
\mathbb V_b=\big\{v(b):v\in\mathbb V\big\},
\end{equation}
and
\begin{equation}\label{eq:Sb}
S_b(v(b),w(b))=-\alpha_v(b)(w(b)),\quad \forall\,v,w\in\mathbb V.
\end{equation}
If $Q$ is contained in $\mathbb V_b$, for instance this is always the case if $t=b$ is not
$(X,\ell_0)$-focal, then we define a symmetric bilinear form
$\mathfrak Q$ on $Q$:
\begin{equation}\label{eq:deffrakQ}
\mathfrak Q = S_Q-S_b\vert_Q.
\end{equation}
\begin{teo}\label{thm:indexthmPQ}
Under the hypotheses of Theorem~\ref{thm:indextheorem}, the index
of $I^\#$ on $\mathcal K^\#$ is finite, and we have:
\begin{equation}\label{eq:indthmPQ}
n_-(I^\#\vert_{\mathcal K^\#})=n_-(B(a)^{-1}\vert_P)+n_-(\mathfrak Q)+
\maslov(X,\ell_0).
\end{equation}
\end{teo}
\begin{proof}
We define the following finite dimensional vector space:
\begin{equation}\label{eq:defVQ}
\mathbb V_Q=\big\{v\in\mathbb V:v(b)\in Q\big\}.
\end{equation}
Since $t=b$ is not focal, it is easily seen that we have a 
direct sum decomposition $\mathcal K^\#=\mathcal K\oplus\mathbb V_Q$. Moreover,
integration by parts in \eqref{eq:defIsust} shows that, for
$v\in \mathbb V_Q$ and $w\in\mathcal H^\#$, it is:
\begin{equation}\label{eq:formII}
I^\#(v,w)=\mathfrak Q(v(b),w(b)).
\end{equation}
This shows that the spaces $\mathcal K$ and $\mathbb V_Q$ are
$I^\#$-orthogonal, and therefore $n_-(I^\#\vert_{\mathcal K^\#})=
n_-(I\vert_{\mathcal K})+n_-(I^\#\vert_{\mathbb V_Q})$.
The conclusion follows from Theorem~\ref{thm:indextheorem}
and from the observation that the isomorphism $\mathbb V_Q\ni v\mapsto v(b)\in Q $
carries $I^\#$ into $\mathfrak Q$.
\end{proof}
\subsection{The opposite symplectic system}
\label{sub:opposite}
It will be useful in the rest of the paper to introduce
an operation of {\em opposition\/} in the category
of symplectic differential systems, that is described briefly in this
subsection.

Let $(X,\ell_0)$ be a set of data for the symplectic differential
problem in $\R^n$; we define a new pair $(X^{\mathrm{op}},\ell_0^{\mathrm{op}})$,
called the {\em opposite symplectic differential problem},
by setting:
\begin{equation}\label{eq:opobjects}
X^\op=\mathcal OX\mathcal O,\quad \ell_0^\op=\mathcal O(\ell_0),
\end{equation}
where $\mathcal O:\R^n\oplus{\R^n}^*\to \R^n\oplus{\R^n}^*$ is given by 
$\mathcal O(v,\alpha)=(v,-\alpha)$.
Opposite symplectic systems have opposite Maslov indexes and opposite index forms; applying
Theorem~\ref{thm:indextheorem} to $(X^\op,\ell_0^\op)$ we get:
\begin{teo}\label{thm:indextheoremop}
Let $(X,\ell_0)$ be a set of data for the symplectic differential problem
in $\R^n$, with $A,B$ of class $C^1$, $C$ continuous, $B(t)$ invertible
for all $t$, and let $\{\mathcal D_t\}_{t\in[a,b]}$ be a $C^2$-family of
$k$-dimensional subspaces of $\R^n$, where $k=n_+(B(t))$, and
$B(t)^{-1}$ is positive definite on $\mathcal D_t$ for all $t$.
Suppose that Assumption~\ref{thm:nondeg} and Assumption~\ref{thm:assK} are satisfied,
and that $t=b$ is not a focal instant. Define  $\mathcal K$ as in \eqref{eq:newK};
then, the co-index of $I$ on $\mathcal K$ is finite, and the following equality
holds:
\begin{equation}\label{eq:equalityidexop}
n_+(I\vert_{\mathcal K})=n_+(B(a)^{-1}\vert_P)-\maslov(X,\ell_0).\qedhere
\end{equation}
\end{teo}

\subsection{Equivalence of Symplectic Differential Systems}
\label{sub:equivalence}
In this subsection we describe the natural isomorphisms in the 
class of symplectic differential problems and we prove the invariance
of the index form by such isomorphisms.

Let $L_0$ be the Lagrangian subspace $\{0\}\oplus{\R^n}^*$ of $(\R^n\oplus{\R^n}^*,
\omega)$; we denote by $\mathrm{Sp}(2n,\R;L_0)$ the closed  subgroup of 
$\mathrm{Sp}(2n,\R)$ consisting of those symplectomorphisms $\phi_0$ such
that $\phi_0(L_0)=L_0$. It is easily seen that any such symplectomorphism
is given in block matrix by:
\begin{equation}\label{eq:phi0}
\phi_0=\left(\begin{array}{cc}Z&0\\ { Z^*}^{-1}W&{
Z^*}^{-1}\end{array}\right),
\end{equation} with $Z:\R^n\to\R^n$ an isomorphism and $W$ a symmetric
bilinear form in $\R^n$.

We give the following:
\begin{defin}\label{thm:defequiv}
The symplectic differential systems with coefficient matrices $X$ and $\tilde X$
 are said to
be {\em isomorphic\/} if there exists a $C^1$-map $\phi_0:[a,b]\to
\mathrm{Sp}(2n,\R;L_0)$ whose upper-left $n\times n$ block is of class $C^2$ and
such that:
\begin{equation}\label{eq:isoX}
\tilde X=\phi_0'\phi_0^{-1}+\phi_0X\phi_0^{-1};
\end{equation}
the pairs  $(X,\ell_0)$ and $(\tilde X,\tilde\ell_0)$  are called
{\em isomorphic\/} if in addition $\tilde\ell_0=\phi_0(a)(\ell_0)$.

We call the map $\phi_0$ an {\em isomorphism\/} between $X$ and
$\tilde X$, or between the pairs $(X,\ell_0)$ and $(\tilde X,\tilde\ell_0)$
in the second case.
\end{defin}
Denoting by $\Psi$ and $\tilde\Psi$ the fundamental matrices of the 
symplectic systems with coefficient matrices $X$ and $\tilde X$ respectively,
we have:
\begin{equation}\label{eq:PsiPsitilde}
\tilde\Psi(t)=\phi_0(t)\Psi(t)\phi_0(a)^{-1},\quad\forall\,t\in[a,b].
\end{equation}
Formula \eqref{eq:PsiPsitilde} is the motivation for Definition~\ref{thm:defequiv}:
the map $\phi_0$ has to be interpreted as a change of variable in the differential
system. Similarly, if $\ell$ and $\tilde\ell$ denote the curve of Lagrangians
associated to the symplectic problems with data $(X,\ell_0)$ and $(\tilde X,\tilde\ell_0)$
respectively, it is easily seen that:
\begin{equation}\label{eq:elltildeell}
\tilde\ell(t)=\phi_0(t)(\ell(t)),\quad\forall\,t\in[a,b].
\end{equation}
Given an isomorphism $\phi_0$ between $(X,\ell_0)$ and $(\tilde X,\tilde\ell_0)$,
we define a $C^2$-curve $Z$ and a $C^1$-curve $W$ by formula \eqref{eq:phi0}. 
The isomorphy between the pairs can be given in terms of the matrix
blocks defining $X$ and $\tilde X$ (recall formula~\eqref{eq:X})
as follows:
\begin{equation}\label{eq:ABC}
\begin{split}
&\tilde A=ZAZ^{-1}-ZBWZ^{-1}+Z'Z^{-1},\\
&\tilde B=ZBZ^*,\\
&\tilde C={Z^*}^{-1}(WA+C-WBW+A^*W+W')Z^{-1}.
\end{split}
\end{equation}
It is easily seen
that $v$ is a $(X,\ell_0)$-solution if and only if $Zv$ is a $(\tilde X,\tilde\ell_0)$
solution. Also, it follows that isomorphic symplectic differential problems
have the same focal instants, and by \eqref{eq:ABC} they have the same signature and
multiplicity. In particular, isomorphic systems have the same focal index.
We also give the formulas for the objects $\tilde P$ and $\tilde S$ defined
by $\tilde\ell_0$ as in \eqref{eq:bijection}:
\begin{equation}\label{eq:PStilde}
\tilde P=Z(a)(P),\quad\tilde
S=S(Z(a)^{-1}\cdot,Z(a)^{-1}\cdot)-W(a)(Z(a)^{-1}\cdot,Z(a)^{-1}\cdot)\vert_{\tilde P}.
\end{equation}
We now prove:
\begin{prop}\label{thm:maslovisom}
Let $(X,\ell_0)$ and $(\tilde X,\tilde\ell_0)$ be isomorphic
pairs such that the final instant is not focal.
Then, $\maslov(X,\ell_0)=\maslov(\tilde X,\tilde\ell_0)$.
\end{prop}
\begin{proof}
Choose $\varepsilon>0$ small enough so that both $\ell$ and
$\tilde\ell$ do not intercept $\Lambda_{\ge1}(L_0)$ on $]a,a+\varepsilon]$.
We prove that $\ell\vert_{[a+\varepsilon,b]}$ and $\tilde\ell\vert_{[a+\varepsilon,b]}$
define the same homology class in $H_1(\Lambda,\Lambda_0(L_0))$. To this aim,
we first observe that the curve $\phi_0$ is homotopic to the constant map
$\phi_0(b)$ in $\mathrm{Sp}(2n,\R;L_0)$. It follows from \eqref{eq:elltildeell}
that $\tilde\ell$ is homotopic to $\phi_0(b)\circ\ell$ in $(\Lambda,\Lambda_0(L_0))$.
It remains to show that the fixed symplectomorphism $\phi_0(b)$ indices the
identity map in the relative homology group $H_1(\Lambda,\Lambda_0(L_0))$.
This follows easily from the fact that, since $\mathrm{Sp}(2n,\R)$ is connected%
\footnote{Observe that the group $\mathrm{Sp}(2n,\R;L_0)$ is {\em not\/} connected,
hence the homotopy argument cannot be done directly at the relative homology level.},
then $\phi_0(b)$ induces the identity map in the absolute homology group
$H_1(\Lambda)$. 
By functoriality of singular homology, it follows that $\phi_0(b)$ also
induces the identity map in the relative homology group $H_1(\Lambda,\Lambda_0(L_0))$.
\end{proof}
Finally, we state the invariance of the index form and of the space $\mathcal K$.
\begin{prop}\label{thm:indexformisom}
Let $(X,\ell_0)$ and $(\tilde X,\tilde\ell_0)$ be isomorphic
pairs, with associated index forms $I:\mathcal H\times
\mathcal H\to\R$ and $\tilde I:\tilde{\mathcal H}\times \tilde{\mathcal H}\to\R$
respectively. Then, the operator $\mathcal H\ni v\mapsto Zv\in\tilde{\mathcal H}$
is a Hilbert space isomorphism that carries $I$ into $\tilde I$ and $\mathcal K$ onto
$\tilde{\mathcal K}$.
\end{prop}
\begin{proof}
It is an easy computation that uses \eqref{eq:ABC} and \eqref{eq:PStilde}.
\end{proof}
Similarly, by a direct computation it is easy to prove the following:
\begin{prop}\label{thm:reducedisomorphic}
Let $X$ be the coefficient matrix of a symplectic differential system in $\R^n$,
and let $\{\mathcal D_t\}_{t\in[a,b]}$ be a $C^2$-family of $k$-dimensional
subspaces of $\R^n$ such that $B(t)$ is negative definite on $\mathcal D_t$ for
all $t$, and $k=n_-(B)$. Then, any two reduced symplectic systems associated
to different choices of bases for $\mathcal D$ are isomorphic.\qed
\end{prop}
\begin{rem}\label{thm:remfacilitar}
In order to ease the verification of Assumption~\ref{thm:assK}, one could study
alternative reduced systems which are isomorphic to \eqref{eq:symplassoc}. An
example of such alternative is given by the following:
\begin{equation}\label{eq:alternativesys}
\left\{\begin{array}{l}f'=-\mathcal B^{-1}\mathcal C\,f+\mathcal
B^{-1}\,\varphi,\\ \varphi'=(\mathcal I-\mathcal C^*\mathcal B^{-1}\mathcal
C)\,f
+\mathcal C^*\mathcal B^{-1}\,\varphi.\end{array}\right.
\end{equation}
A sufficient condition for a symplectic differential problem not have focal instants
is given in Remark~\ref{thm:remposdef}. We remark that this condition is {\em not\/}
preserved by isomorphisms; hence, another sufficient condition
that guarantees the validity of Assumption~\ref{thm:assK} is that
the symmetric bilinear form $\mathcal I-\mathcal C^*\mathcal B^{-1}\mathcal
C$ be negative semi-definite.
\end{rem}

Proposition~\ref{thm:indexformisom} and Proposition~\ref{thm:reducedisomorphic}
show that the index form  of any reduced symplectic system defines
a unique symmetric bilinear form on the space $\mathcal S$ of $\mathcal D$-valued
vector fields vanishing at the endpoints.  
A simple computation shows that this reduced index form
is the restriction of the original index form:
\begin{prop}\label{thm:restrindexform}
Let $(X,\ell_0)$ be a set of data for the symplectic differential
problem  in $\R^n$, $\mathcal D$ be a $C^2$-family of subspaces
of $\R^n$, and let $Y_1,\ldots,Y_k$ be a basis of $\mathcal D$. Denote
by $I$ the index form of $(X,\ell_0)$ and
by $I_{\mathrm{red}}$ the index form of the reduced symplectic system of $(X,\ell_0)$
associated to the given basis. Then, for every $f,g\in H^1_0([a,b];\R^k)$,
we have:
\begin{equation}\label{eq:restrindexform}
I_{\mathrm{red}}(f,g)=I(v,w),
\end{equation}
where $v=\sum_if_iY_i$ and $w=\sum_ig_iY_i$.\qed
\end{prop}
\begin{cor}\label{thm:corIdefnegS}
Assumption~\ref{thm:assK} is equivalent to the condition that the 
restriction of the index form
$I$ to the space $\mathcal S$ (defined in \eqref{eq:defS}) be negative definite.
\end{cor}
\begin{proof}
Apply Theorem~\ref{thm:indextheoremop} to the reduced
symplectic differential system, keeping in mind Proposition~\ref{thm:restrindexform}.
\end{proof}
\subsection{Every symplectic system is isomorphic to a Morse--Sturm
system}
\label{sub:isoMS}
It is a trivial observation that every second order linear
differential equation $x''+f_1 x'+f_2 x=0$ can be written as
a Sturm equation:   $(e^{\int f_1} x')'=-f_2e^{\int f_1}x$. 
On the other hand, given any Sturm equation $(px')'=rx$, with $p$ a map
of class $C^2$ and, say,
$p>0$, then  the change of variable $y=\sqrt px$ transforms the
equation into the Sturm equation $y''=\tilde ry$, with
$\tilde r=\left[\frac rp+\frac{p''p-\frac12 (p')^2}{2p^2}\right]$.

In the language of symplectic systems, these facts are expressed
by saying that every  $C^2$
unidimensional symplectic system is isomorphic
to a Morse--Sturm system having coefficient matrix $\left(\begin{array}{cc}
0&b\\ c&0\end{array}\right)$ with $b$ constant. 

This fact holds for symplectic systems of any dimension:
\begin{prop}\label{thm:isoMS}
Every symplectic differential system \eqref{eq:sympl} in $\R^n$  
is isomorphic to a
Morse--Sturm system, i.e., a symplectic system with coefficient matrix $\tilde
X=\left(\begin{array}{cc}
0&\tilde B\\ \tilde C&0\end{array}\right)$. Moreover, if $B$ is a map
of class $C^2$, then $\tilde B$ can be taken to be the {\em constant\/}
map\begin{equation}\label{eq:constantmap}
\left(\begin{array}{cc}-\mathrm{Id}_k&0\\0&\mathrm{Id}_{n-k}\end{array}\right), 
\end{equation}
where $k=n_-(B)$.
\end{prop}
\begin{proof}
We will use the notation of Subsection~\ref{sub:equivalence} to exhibit the required
isomorphisms. For the first part of the statement, simply consider the symmetric matrix $W=0$
and let $Z$ be the solution on $[a,b]$ of the initial value problem $Z'=-ZA$ with
$Z(a)=\mathrm{Id}$. Then, by formula \eqref{eq:ABC}, we have $\tilde A=0$.

For the second part of the thesis, observe first that if $B$ is a map of class
$C^2$, then we can assume that $B$ is the constant map \eqref{eq:constantmap},
observing that we can always find a basis $\{e_1(t),\ldots,e_n(t)\}$
in $\R^n$ such that $B(t)$ is in the canonical diagonal form \eqref{eq:constantmap}
and such that each $e_i(t)$ is of class $C^2$.

Denote by $\mathrm O(k,n-k)$ the closed subgroup of $\mathrm{GL}(n,\R)$
consisting of those maps preserving the quadratic form with matrix
\eqref{eq:constantmap}  and by $\mathrm o(k,n-k)$ its Lie algebra. 
Then, if $\phi_0$ is an
isomorphism between the symplectic systems $X$ and
$\tilde X$ such that the corresponding map $Z$ takes 
values in $\mathrm O(k,n-k)$, it follows that
$\tilde B=ZBZ^*=B$ is also equal to \eqref{eq:constantmap}. Now, if we consider the
symmetric matrix
$W=\frac12(BA+A^*B)$ and if we take $Z$ to be the solution of the initial value problem
$Z'=Z(BW-A)$, $Z(a)=\mathrm{Id}$, then $\tilde A=0$, and, since
$BW-A\in\mathrm o(k,n-k)$, we obtain $Z\in\mathrm O(k,n-k)$. This concludes the proof.
\end{proof}
Examples of Morse--Sturm systems in $\R^n$ whose coefficient matrix has constant
upper-right $n\times n$ block arise from Jacobi equations along semi--Riemannian
geodesics by a parallel trivialization of the tangent bundle
along the geodesic (see Subsection~\ref{sub:semiRiemannian}). 
Conversely, every Morse--Sturm system of
this type corresponds to the Jacobi equation along a semi--Riemannian geodesic:
\begin{prop}\label{thm:MSJacobi}
Every smooth Morse--Sturm system in $\R^n$ whose coefficient matrix
has  upper-right $n\times n$ block constant can be obtained 
as the Jacobi equation along a geodesic $\gamma$ in a semi-Riemannian
manifold $(M,g)$ by a parallel trivialization of the tangent bundle $TM$
along $\gamma$. 
\end{prop}
\begin{proof}
See \cite[Section~3]{Hel1} and also \cite[Proposition~2.3.1]{MPT}.
\end{proof}
\end{section}


\begin{section}{Applications to Hamiltonian Systems}
\label{sec:hamiltonian}
In this section we will consider the following setup.
Let $(\mathcal M,\omega)$ be a symplectic manifold, i.e., $\mathcal M$
is a smooth manifold and $\omega$ is a smooth closed skew-symmetric
nondegenerate two-form on $\mathcal M$; we set $\mathrm{dim}(\mathcal M)=2n$.
Let $H:U\to\R$ be a smooth function defined in an open set $U\subseteq\R\times\mathcal
M$; we will call such function a {\em Hamiltonian\/} in $(\mathcal M,\omega)$.
For each $t\in\R$, we denote by $H_t$ the map $m\mapsto H(t,m)$ defined in
the open set $U_t\subseteq\mathcal M$ consisting of those $m\in\mathcal M$ such that
$(t,m)\in U$.
We denote by $\Hf$ the smooth time-dependent vector field in  
$\mathcal M$ defined by $\mathrm dH_t(m)=\omega(\Hf(t,m),\cdot)$ for all
$(t,m)\in U$; let $F$ denote the maximal flow of the vector field $\Hf$
defined on an open set of $\R\times\R\times\mathcal M$ taking values in
$\mathcal M$, i.e., for each $m\in\mathcal M$ and $t_0\in\R$, 
the curve $t\mapsto F(t,t_0,m)$ is
a maximal integral curve of $\Hf$ and $F(t_0,t_0,m)=m$. This means that $F(\cdot,t_0,m)$
is a maximal solution of the equation:
\[\frac{\mathrm d}{\mathrm dt}\,F(t,t_0,m)=\Hf(t,F(t,t_0,m)),\quad F(t_0,t_0,m)=m.\]
Recall that $F$ is a smooth map; we also write $F_{t,t_0}$ for the map $m\mapsto F(t,t_0,m)$;
observe that $F_{t,t_0}$ is a diffeomorphism between open subsets of $\mathcal M$.

We recall that a {\em symplectic chart\/} in $\mathcal M$ is a local chart $(q,p)$
taking values in $\R^n\oplus{\R^n}^*$ whose differential at each point
is a symplectomorphism from the tangent space of $\mathcal M$ to $\R^n\oplus{\R^n}^*$
endowed with the canonical symplectic structure. We write $q=(q_1,\ldots,q_n)$
and $p=(p_1,\ldots,p_n)$; we denote by $\{\frac{\partial}{\partial
q_i},\frac{\partial}{\partial p_j}\}$, $i,j=1,\ldots,n$ the corresponding local referential
of $T\mathcal M$, and by $\{\mathrm dq_i,\mathrm dp_j\}$ the local referential
of $T\mathcal M^*$. By Darboux's Theorem, there always
exists an atlas of symplectic charts.

In a given symplectic chart $(q,p)$, we have:
\[\omega=\sum_{i=1}^n\mathrm dq_i\wedge\mathrm dp_i,\quad\Hf=\sum_{i=1}^n\left(\frac{\partial
H}{\partial p_i}\,\frac\partial{\partial q_i}-\frac{\partial
H}{\partial q_i}\,\frac\partial{\partial p_i}\right).\]

Let $\mathcal P$ be a {\em Lagrangian submanifold\/} of $\mathcal M$, i.e.,
$T_m\mathcal P$ is a Lagrangian subspace of $T_m\mathcal M$ for every $m\in \mathcal P$.
We fix an integral curve $\Gamma:[a,b]\to\mathcal M$ of $\Hf$, so that
$\Gamma(t)=F(t,a,\Gamma(a))$ for all $t\in[a,b]$. 
We also say that $\Gamma$ is a {\em solution of
the Hamilton  equations}, i.e., in a symplectic chart $\Gamma(t)=(q(t),p(t))$:
\begin{equation}\label{eq:HamJac}
 \frac{\mathrm
dq}{\mathrm dt}= \frac{\partial H}{\partial p}, \quad
\frac{\mathrm dp}{\mathrm dt}=-\frac{\partial H}{\partial q}.
\end{equation}
We assume that $\Gamma$ starts at $\mathcal P$, that is $\Gamma(a)\in\mathcal P$.
Finally, we will consider a fixed smooth distribution
$\mathfrak L$ in $\mathcal M$ such that $\mathfrak L_m$ is a Lagrangian
subspace of $T_m\mathcal M$ for all $m\in\mathcal M$.

The basic example to keep in mind for the above setup is the case
where $\mathcal M$ is the cotangent bundle $TM^*$ of some smooth
manifold $M$ endowed with the canonical symplectic structure, 
$\mathcal P$ is the annihilator $TP^o$ of some smooth submanifold
$P$ of $M$, and $\mathfrak L$ is the distribution consisting of the
{\em vertical\/} subspaces, i.e., the subspaces tangent to the fibers
of $TM^*$. All the results proven in the abstract framework will
be discussed in detail for this special case in Subsection~\ref{sub:hyper}.

\subsection{The symplectic differential problem associated to a Hamiltonian
system}\label{sub:associated}
It is well known that the Hamiltonian flow $F_{t,t_0}$ consists of
symplectomorphisms:
\begin{prop}\label{thm:omegainvariant}
The symplectic form $\omega$ is invariant by the Hamiltonian
flow $F$, i.e., $F_{t,t_0}^*\omega=\omega$ for all $(t,t_0)$.\qed
\end{prop}
A sextuplet $\sextuple$
where $(\mathcal M,\omega)$ is a symplectic manifold, $H$ is a (time-dependent)
Hamiltonian function defined on an open subset of $\R\times\mathcal M$,
$\mathfrak L$ is a smooth distribution of Lagrangians in $\mathcal M$,
$\Gamma:[a,b]\to\mathcal M$ is an integral curve of $\Hf$ and $\mathcal P$
is a Lagrangian submanifold of $\mathcal M$ with $\Gamma(a)\in\mathcal P$,
will be called {\em a set of data for the Hamiltonian problem}.

We   give some more basic definitions.

\begin{defin}\label{thm:deflinsol}
A vector field $\rho$ along $\Gamma$ in $\mathcal M$ is said to be a {\em solution 
for the linearized Hamilton  (LinH) equations\/} if it satisfies:
\begin{equation}\label{eq:linsol}
\rho(t)=\mathrm dF_{t,a}(\Gamma(a))\,\rho(a).
\end{equation}
We also say that $\rho$ is a {\em $\mathcal P$-solution\/}  for the (LinH) equations if in
addition it satisfies $\rho(a)\in T_{\Gamma(a)}\mathcal P$.
\end{defin}
The solutions of the (LinH) equations are precisely
the variational vector fields along $\Gamma$ corresponding to variations
of $\Gamma$ by integral curves of $\Hf$; the $\mathcal P$-solutions correspond 
to variations by integral curves starting on $\mathcal P$.
It is easy to see that the set of solutions of the (LinH) equations
is a vector space of dimension $2n$ and that the $\mathcal P$-solutions
form an $n$-dimensional vector subspace.

\begin{defin}\label{thm:defPfocal}
A point $\Gamma(t)$, $t\in]a,b]$  is said to be a {\em $\mathcal P$-focal
point\/} along $\Gamma$ if there exists a non zero $\mathcal P$-solution $\rho$ for the 
(LinH) equations such that $\rho(t)\in\mathfrak L_{\Gamma(t)}$.
The {\em multiplicity\/} of a $\mathcal P$-focal point $\Gamma(t)$ is the dimension
of the vector space of such $\rho$'s.
\end{defin}
We will now describe how, under suitable nondegeneracy hypotheses, one associates 
to the set of data for the Hamiltonian problem $\sextuple$ 
an isomorphism class of sets of data for the
symplectic differential problem.
\begin{defin}\label{thm:defLsympl}
A {\em symplectic $\mathfrak L$-trivialization\/} of $T\mathcal M$ along $\Gamma$
is a smooth family of symplectomorphisms $\phi(t):\R^n\oplus{\R^n}^*
\to T_{\Gamma(t)}\mathcal M$ such that, for all $t\in[a,b]$, $\phi(t)(L_0)=\mathfrak
L_{\Gamma(t)}$, where $L_0=\{0\}\oplus{\R^n}^*$.
\end{defin}
The existence of symplectic $\mathfrak L$-trivializations along $\Gamma$ is easily
established with elementary arguments, using the fact that $T\mathcal M$
restricts to a trivial vector bundle along $\Gamma$. 

We will be interested also in the quotient bundle $T\mathcal M/\mathfrak L$
and its dual bundle. We have an obvious canonical identification of
the dual $(T\mathcal M/\mathfrak L)^*$ with the annihilator
$\mathfrak L^o\subset T\mathcal M^*$; moreover, using the symplectic form,
we will identify $\mathfrak L^o$ with $\mathfrak L$ by the isomorphism:
\begin{equation}\label{eq:idTMTM*}
T_m\mathcal M\ni\rho\mapsto \omega(\cdot,\rho)\in T_m\mathcal M^*,\quad
m\in\mathcal M.
\end{equation} 
A symplectic $\mathfrak L$-trivialization $\phi$ induces a trivialization
of the quotient bundle $T\mathcal M/\mathfrak L$ along $\Gamma$, namely,
for each $t\in[a,b]$ we define an isomorphism $\mathcal Z_t:\R^n\to
T_{\Gamma(t)}\mathcal M /\mathfrak L_{\Gamma(t)}$:
\begin{equation}\label{eq:Zcal}
\mathcal Z_t(x)=\phi(t)(x,0)+\mathfrak L_{\Gamma(t)},\quad x\in\R^n.
\end{equation}

Given a
symplectic $\mathfrak L$-trivialization $\phi$ of $T\mathcal M$ along $\Gamma$,  
we define a smooth curve $\Psi:[a,b]\to\mathrm{Sp}(2n,\R)$ by:
\begin{equation}\label{eq:Psitriv}
\Psi(t)=\phi(t)^{-1}\circ\mathrm dF_{t,a}(\Gamma(a))\circ\phi(a).
\end{equation}
The fact that $\Psi(t)$ is a symplectomorphism follows from
Proposition~\ref{thm:omegainvariant}.

We now define a smooth curve $X:[a,b]\to\mathrm{sp}(2n,\R)$ by setting:
\begin{equation}\label{eq:Xtriv}
X(t)=\Psi'(t)\Psi(t)^{-1};
\end{equation}
The $n\times n$ blocks of the matrix $X$ will be denoted by $A,B$ and $C$, as in
formula \eqref{eq:X}. Finally, we define a Lagrangian subspace
$\ell_0$ of $\R^n\oplus{\R^n}^*$ by:
\begin{equation}\label{eq:ell0triv}
\ell_0=\phi(a)^{-1}(T_{\Gamma(a)}\mathcal P).
\end{equation}
Suppose now that we are given another symplectic $\mathfrak L$-trivialization
$\tilde\phi$ of $T\mathcal M$ along $\gamma$. Denote by $\tilde{\mathcal Z}$ the
relative trivialization of $T\mathcal M/\mathfrak L$, and by $\tilde\Psi,
\tilde X$ and $\tilde\ell_0$ the objects defined in \eqref{eq:Psitriv},
\eqref{eq:Xtriv} and \eqref{eq:ell0triv} relatively to $\tilde\phi$. 

If we define the smooth curve $\phi_0:[a,b]\to\mathrm{Sp}(2n,\R;L_0)$
by:
\[\phi_0(t)=\tilde\phi(t)^{-1}\circ\phi(t),\]
then obviously $\Psi$ and $\tilde\Psi$ are related as in 
formula \eqref{eq:PsiPsitilde}, and $\ell_0$ is related with $\tilde\ell_0$
as in \eqref{eq:elltildeell}. Differentiating \eqref{eq:PsiPsitilde}, we obtain
easily \eqref{eq:isoX}. Moreover, if we set:
\[Z_t=\tilde{\mathcal Z}_t^{-1}\circ\mathcal Z_t,\]
then, $Z(t)$ is the upper-left $n\times n$ block of $\phi_0(t)$, as in
\eqref{eq:phi0}. 

We are ready for the following:
\begin{defin}\label{thm:defFCH}
The {\em canonical bilinear form\/}  of the set of data 
$\sextuple$ is a family of symmetric bilinear forms
$\HL(t)$ on $(T_{\Gamma(t)}\mathcal M/\mathfrak L_{\Gamma(t)})^*\simeq \mathfrak L_{\Gamma(t)}^o
\simeq  \mathfrak L_{\Gamma(t)}$ 
given by:
\begin{equation}\label{eq:defHL}
\HL(t)=\mathcal Z_t\circ B(t)\circ\mathcal Z_t^*,
\end{equation}
where $\mathcal Z$ is the trivialization of $T\mathcal M/\mathfrak L$ relative
to some symplectic $\mathfrak L$-trivialization $\phi$ of
$T\mathcal M$ and $B$ is the upper-right $n\times n$ block of
the map $X$ in \eqref{eq:Xtriv}. Observe that, by \eqref{eq:ABC}
and the construction of the map $\mathcal Z$, the right hand side
of \eqref{eq:defHL} does {\em not\/} depend on the choice of
the symplectic $\mathfrak L$-trivialization of $T\mathcal M$.

We say that the set of data $\sextuple$
is {\em nondegenerate\/} if $\HL(t)$ is nondegenerate for all $t\in[a,b]$.
In this case, we can also define the symmetric bilinear form $\HL(t)^{-1}$ on
$ T_{\Gamma(t)}\mathcal M/\mathfrak L_{\Gamma(t)}$.
\end{defin}
Given a nondegenerate set of data $\sextuple$, let us consider the pair $(X,\ell_0)$ defined
by \eqref{eq:Xtriv} and \eqref{eq:ell0triv}.
It is easily seen that the submanifold $\mathcal P$ and the space
$P$ defined by $\ell_0$ as in \eqref{eq:bijection} are related by the following:
\[\mathcal Z_a(P)=\pi(T_{\Gamma(a)}\mathcal P),\]
where $\pi:T_{\Gamma(a)}\mathcal M\to
T_{\Gamma(a)}\mathcal M/\mathfrak L_{\Gamma(a)}$ is the quotient map.
We set:
\begin{equation}\label{eq:Pcal0}
\mathcal P_0=\pi(T_{\Gamma(a)}\mathcal P).
\end{equation}
In analogy with Assumption~\ref{thm:nondeg}, we make the following
\begin{assum}\label{thm:nondegH}
We assume that the symmetric bilinear form $\HL(a)^{-1}$ is nondegenerate
on $\mathcal P_0$. This is equivalent to requiring
that $\HL(a)$ is nondegenerate on the annihilator $(T_{\Gamma(a)}\mathcal P+
\mathfrak L_a)^o\subset\mathfrak L_a^o\subset T_{\Gamma(a)}\mathcal M^*$.
\end{assum}
Given a nondegenerate set of data $\sextuple$ satisfying
Assumption~\ref{thm:nondegH}, the pair $(X,\ell_0)$ defined by  
\eqref{eq:Xtriv} and \eqref{eq:ell0triv} is a set of data for the symplectic
differential problem. We say that $(X,\ell_0)$ is the pair {\em associated\/}
to $\sextuple$ by the choice of the symplectic $\mathfrak L$-trivialization
$\phi$ of $T\mathcal M$ along $\Gamma$. We have proven above that 
pairs $(X,\ell_0)$ and $(\tilde X,\tilde\ell_0)$ associated to $\sextuple$ by
different choices of symplectic $\mathfrak L$-trivializations are isomorphic.

To define the signature of a $\mathcal P$-focal point along $\Gamma$, we need
to introduce the following space, in analogy with formula \eqref{eq:spaces}:
\begin{equation}\label{eq:deffrakA}
\mathfrak V[t]=\big\{\rho(t):\rho\ \text{is a $\mathcal P$-solution of the
(LinH) equation}\big\}\cap\mathfrak L_{\Gamma(t)},\quad
t\in[a,b].
\end{equation}
Using the isomorphism \eqref{eq:idTMTM*}, it is easy to see that $\mathfrak V[a]$
is identified with the annihilator $(T_{\Gamma(a)}\mathcal P+
\mathfrak L_a)^o$. It is easily seen that a point $\Gamma(t)$ is $\mathcal P$-focal
if and only if $\mathfrak V[t]$ is not zero and that the dimension of
$\mathfrak V[t]$ is precisely the multiplicity of $\Gamma(t)$.
\begin{defin}\label{thm:defsgnPfocal}
Let $\Gamma(t)$ be a $\mathcal P$-focal point along $\Gamma$. The 
{\em signature\/} $\sgn(\Gamma(t))$ is the signature of the restriction
of $\HL(t)$ to $\mathfrak V[t]\subset\mathfrak L_t\simeq\mathfrak L_t^o$.
$\Gamma(t)$ is said to be a nondegenerate $\mathcal P$-focal point if 
such restriction is nondegenerate. If $\Gamma$ has only a finite
number of $\mathcal P$-focal points, we define the {\em focal index\/}
$\mathrm i_{\mathrm{foc}}(\Gamma)$ as:
\begin{equation}\label{eq:deffocalindGamma}
\mathrm i_{\mathrm{foc}}(\Gamma)=\sum_{t\in\,]a,b]}\sgn(\Gamma(t)).
\end{equation}
\end{defin}
\subsection{The index theorem for Hamiltonian systems}
\label{sub:maslovH}
Using the identification with isomorphism
classes of symplectic problems, we now translate the result of the
Index Theorem~\ref{thm:indextheorem} for nondegenerate sextuplets
$\sextuple$.
The sextuplets considered henceforth will be assumed to be nondegenerate
and satisfying Assumption~\ref{thm:nondegH}.

Recalling that the Maslov index of isomorphic symplectic problems
are equal (Proposition~\ref{thm:maslovisom}), we can give the following:
\begin{defin}\label{thm:defMaslovH}
Given a  set of data $\sextuple$ such that $\Gamma(b)$ is not a $\mathcal P$-focal point, we define its
{\em Maslov index\/} $\maslov(\Gamma)$ as the Maslov index of any pair
$(X,\ell_0)$ associated to it by a symplectic $\mathfrak L$-trivialization
of $T\mathcal M$ along $\Gamma$.
\end{defin}
Recalling Proposition~\ref{thm:indexformisom}, we can also define an index
form associated to a Hamiltonian problem as follows.
\begin{defin}\label{thm:defindexformH}
Let $\mathcal H_{\Gamma}$ denote the Hilbert space of sections
$\mathfrak v$ of Sobolev class $H^1$ of  the quotient bundle
$T\mathcal M/\mathfrak L\/$ along $\Gamma$ such that $\mathfrak v(a)\in\mathcal P_0$ and
$\mathfrak v(b)=0$. The {\em index form\/}
$I_\Gamma$ of the Hamiltonian problem is the bounded symmetric
bilinear form on $\mathcal H_{\Gamma}$ defined by:
\begin{equation}\label{eq:defindexformH}
I_\Gamma(\mathfrak v,\mathfrak w)=I(\mathcal Z^{-1}\mathfrak v,\mathcal Z^{-1}\mathfrak
w),
\end{equation}
where $\mathcal Z$ is the trivialization of $T\mathcal M/\mathfrak L$
along $\Gamma$ associated to a symplectic $\mathfrak L$-trivialization
$\phi$ of $T\mathcal M$ as in \eqref{eq:Zcal}, and $I$ is the index
form of the pair $(X,\ell_0)$ associated to such trivialization.
\end{defin}

Setting $k=n_-(\HL)$, we will now consider a smooth family of $k$-dimensional 
subspaces $\mathcal D_t\subset T_{\Gamma(t)}\mathcal M/\mathfrak L_{\Gamma(t)}$,
$t\in[a,b]$ (this is a subbundle of the pull-back
$\Gamma^*(T\mathcal M/\mathfrak L)$). We assume that $\HL(t)^{-1}$ is
negative definite on $\mathcal D_t$ for all $t$. We call such a family a
{\em maximal negative distribution along $\Gamma$}. Given the distribution
$\mathcal D$, we define the closed subspaces $\mathcal S_{\mathcal D},\mathcal
K_{\mathcal D}\subset
\mathcal H_{\Gamma}$ by:
\begin{equation}\label{eq:SD}
\mathcal S_{\mathcal D}=\big\{\mathfrak v\in \mathcal H_{\Gamma}:
\mathfrak v(a)=\mathfrak v(b)=0,\ \mathfrak v(t)\in\mathcal D_t,\ \forall\,t\big\},
\end{equation}
\begin{equation}\label{eq:KD}
\mathcal K_{\mathcal D}=\big\{\mathfrak v\in \mathcal H_{\Gamma}:
\mathcal Z^{-1}(\mathfrak v)\in \mathcal K\big\},
\end{equation}
where $\mathcal Z$ is the trivialization of $T\mathcal M/\mathfrak L$
along $\Gamma$ associated to a symplectic $\mathfrak L$-trivialization
$\phi$ of $T\mathcal M$ as in \eqref{eq:Zcal}, and $\mathcal K$ is defined
as in \eqref{eq:newK} corresponding to the distribution $\mathcal Z^{-1}(\mathcal D)$
of subspaces of $\R^n$. By Proposition~\ref{thm:indexformisom}, the definition
of the space $\mathcal K_{\mathcal D}$ does not depend on the choice
of the trivialization $\phi$.

\begin{teo}\label{thm:main}
Let $\sextuple$ be a nondegenerate set of data
for the Hamiltonian problem satisfying Assumption~\ref{thm:nondegH}
and such that $\Gamma(b)$ is not a $\mathcal P$-focal point. Let $\mathcal D$
be a maximal negative distribution along $\Gamma$, and suppose that the
index form $I_\Gamma$ is {\em negative definite\/} on the space
$\mathcal S_{\mathcal D}$ defined in \eqref{eq:SD}. Let's denote by $\mathcal S_{\mathcal
D}^\perp$ the $I_\Gamma$-orthogonal complement of $\mathcal S_{\mathcal D}$ in
$\mathcal H_{\Gamma}$, i.e., 
\begin{equation}\label{eq:SDperp}
\mathcal S_{\mathcal D}^\perp=\big\{\mathfrak v\in\mathcal H_{\Gamma}:I_\Gamma(\mathfrak
v,\mathfrak w)=0,\ \ \forall\,\mathfrak w\in
\mathcal S_{\mathcal D}\big\}.
\end{equation}
Then, $\mathcal H_{\Gamma}=\mathcal S_{\mathcal D}^\perp\oplus\mathcal S_{\mathcal D}$,
$\mathcal S_{\mathcal D}^\perp=\mathcal K_{\mathcal D}$ and:
\begin{equation}\label{eq:main}
n_-(I_\Gamma\vert_{\mathcal K_{\mathcal D}})=n_-(\HL(a)^{-1}\vert_{
\mathcal P_0})+\maslov(\Gamma).
\end{equation}
\end{teo}
\begin{proof}
Use a symplectic $\mathfrak L$-trivialization of $T\mathcal M$ along $\Gamma$
to associate to the given sextuplet $\sextuple$ a set of data for the symplectic differential
problem, for which the results of Section~\ref{sec:symplectic} apply.
The thesis follows from Corollary~\ref{thm:H=K+S}, Lemma~\ref{thm:SKorth},
Theorem~\ref{thm:indextheorem} and Corollary~\ref{thm:corIdefnegS}.
\end{proof}
\subsection{The case of a hyper-regular Hamiltonian system}
\label{sub:hyper}
We will now consider the special case of a Hamiltonian problem defined
in the cotangent bundle $\mathcal M=TM^*$ of a smooth $n$-dimensional 
manifold $M$ endowed
with the canonical symplectic form $\omega$. Recall that $\omega$ is
defined by $\omega=-\mathrm
d\theta$, where $\theta$ is the canonical $1$-form on $TM^*$, given by $\theta_p(\rho)=
p(\mathrm d\pi_p(\rho))$ where $\pi:TM^*\to M$ is the projection,
$p\in TM^*$, $\rho\in T_pTM^*$.

We will prove that, in this situation, if the Hamiltonian is hyper-regular,
i.e., if it corresponds to a Lagrangian, then the index form $I_\Gamma$ is
the second variation of the Lagrangian action principle.
In the general context of subsection~\ref{sub:associated} we do
not know whether there exists an interpretation of $I_\Gamma$ as
a second variation of some functional.

A local chart $(q_1,\ldots,q_n)$ on $M$ induces a local chart
$(q_1,\ldots,q_n,p_1,\ldots,p_n)$ on $TM^*$, which is symplectic.
We consider the Lagrangian distribution $\mathfrak L$ on $TM^*$
given by:
\begin{equation}\label{eq:verticaltang}
\mathfrak L_p=\mathrm{Ker}(\mathrm
d\pi_p)=T_p(T_{\pi(p)}M^*)\simeq T_{\pi(p)}M^*,\quad p\in TM^*.
\end{equation}
We also call $\mathfrak L_p $ the {\em vertical subspace\/} of $T_pTM^*$.

Let $P$ be a smooth submanifold of $M$; then the annihilator
\[\mathcal P=TP^o=\big\{p\in TM^*:\pi(p)=m\in P\ \text{and}\ p\vert_{T_{m}P}=0\big\}\]
is a Lagrangian submanifold of $TM^*$. Indeed, $\theta$ vanishes on $TP^o$
and so does $\omega$. 

We consider a Hamiltonian $H$ defined on an open subset of $\R\times TM^*$
and an integral curve $\Gamma:[a,b]\to TM^*$ of $\Hf$ such that
$\Gamma(a)\in TP^o$; let $\gamma$ be the projection in $M$ of $\Gamma$:
\[\gamma=\pi\circ\Gamma,\]
so that $\gamma(a)\in P$.

Using the projection $\pi$, we can identify the quotient $T_p\mathcal M/\mathfrak L_p$
with the tangent space $T_{\pi(p)}M$, and the space $\mathcal P_0$ defined
in \eqref{eq:Pcal0} becomes identified with $T_{\gamma(a)}P$. By these identifications,
we can describe the Hilbert space $\mathcal H_{\Gamma}$ of
Definition~\ref{thm:defindexformH}  as the space of vector fields
$\mathfrak v$ along $\gamma$ in $M$, of Sobolev class $H^1$, such that 
$\mathfrak v(a)\in T_{\gamma(a)}P$ and $\mathfrak v(b)=0$.

Let $q=(q_1,\ldots,q_n)$ be a local chart in $M$ and let
$(q,p)=(q_1,\ldots,q_n,p_1,\ldots,p_n)$ be the corresponding chart in $TM^*$; suppose that we
have a subinterval $[c,d]\subset[a,b]$ such that $\gamma([c,d])$ is contained in the domain
of the chart $q$. For each $t\in[c,d]$, the differential of the chart $(q,p)$ gives
a symplectomorphism from $\R^n\oplus{\R^n}^*$ to $T_{\Gamma(t)}TM^*$ which takes
$L_0=\{0\}\oplus{\R^n}^*$ to the vertical space $\mathfrak L_{\Gamma(t)}$.
Suppose that $\phi$ is a symplectic $\mathfrak L$-trivialization of $TTM^*$
along $\Gamma$ such that, for $t\in[c,d]$, $\phi(t)$ is the symplectomorphism
induced by the chart $(q,p)$. In this case we say that the
symplectic $\mathfrak L$-trivialization is {\em compatible\/} with the
chart $(q,p)$ on the interval $[c,d]$; the corresponding trivialization $\mathcal Z$
of the quotient $T\mathcal M/\mathfrak L$ in the interval $[c,d]$ is given by
the differential of the chart $q$. 
Now, consider the pair
$(X,\ell_0)$ of data for the symplectic differential problem associated to $\phi$; for
$t\in[c,d]$, we have:
\begin{equation}\label{eq:XH}
X=\begin{pmatrix}\frac{\partial^2 H}{\partial
q\partial p}& \frac{\partial^2 H}{\partial p^2}\\  
-\frac{\partial^2 H}{\partial q^2}&  -\frac{\partial^2 H}{\partial
p\partial q}\end{pmatrix}.
\end{equation}
This is obtained by observing that solutions $(v,\alpha)$ of the symplectic 
differential system with coefficient matrix \eqref{eq:XH} correspond by $\phi$ to
solutions $\rho$  of the linearized Hamilton  equations; such $\rho$'s 
are variational vector fields corresponding to variations of $\Gamma$ by 
integral curves of $\Hf$.
In this setup, the canonical bilinear form $\HL(t)$ of \eqref{eq:defHL} is
a bilinear form in $T_{\gamma(t)}M^*$, and we have the following:
\begin{prop}\label{thm:HLTM*}
The canonical bilinear form $\HL(t)$ equals the Hessian at $\Gamma(t)$ of the map
$p\mapsto H(t,p) $, defined in a neighborhood of $\Gamma(t)$ in the fiber
$T_{\gamma(t)}M^*$.
\end{prop}
\begin{proof}
It follows directly from \eqref{eq:XH}.
\end{proof}
Let us now assume that $H$ is a {\em hyper-regular\/} Hamiltonian on an open
subset $U\subset\R\times TM^*$, i.e.,
its fiber derivative is a (global) diffeomorphism between $U$ and an open
subset $V\subset\R\times TM$ (see \cite{AM}); let $L$ be the corresponding
hyper-regular Lagrangian on $V$.
The map $L$ defines an {\em action functional\/}
on the set of $C^1$-curves  $\vartheta:[a,b]\to M$ such that $(t,\vartheta'(t))\in
V$ for all $t\in[a,b]$, by:
\begin{equation}\label{eq:actionfunct}
\vartheta\longmapsto\int_a^bL(t,\vartheta'(t))\;\mathrm dt.
\end{equation}
The set of such curves $\vartheta$ has the structure of an infinite dimensional
Banach manifold, and the action functional above is smooth.
It is well known (see \cite[Chapter~3]{AM}) that the critical points of the restriction
of the action functional \eqref{eq:actionfunct} to the space of $C^1$-curves
connecting the submanifold $P$ and the point $\gamma(b)$ are the curves
$\vartheta:[a,b]\to M$ such that $\vartheta(b)=\gamma(b)$ and  $t\mapsto 
\mathrm d\big(L\vert_{(\{t\}\times T_{\vartheta(t)}M)\cap V}\big)
(t,\vartheta'(t))$ is an integral curve of $\Hf$ starting at $T\mathcal P_0$.

Hence, $\gamma$ is a critical point of the restriction of the action functional,
and the Hessian of such functional at $\gamma$ defines a bounded symmetric
bilinear form on the Banach space of $C^1$-vector fields $\mathfrak v$ along
$\gamma$ such that $\mathfrak v(a)\in TP$ and $\mathfrak v(b)=0$. This Hessian
is computed in the following:
\begin{prop}\label{thm:Hessianaction}
The Hessian of the action functional \eqref{eq:actionfunct} at $\gamma$ 
has a (unique) extension to a  bounded symmetric bilinear form on the Hilbert
space $\mathcal H_{\Gamma}$, which is given by $I_\Gamma$ (see
Definition~\ref{thm:defindexformH}).
\end{prop}
\begin{proof}
If $\mathfrak v,\mathfrak w$ are vector fields along $\gamma$ with
``small'' support, then the computation of $I_\Gamma(\mathfrak v,\mathfrak w)$
can be done in local coordinates using \eqref{eq:XH}. In this case,
the conclusion is easily obtained.
For the general case, one simply observes that both the index form and the Hessian
are bilinear, and that any smooth vector field along $\gamma$ can be written as
a finite sum of vector fields with ``small'' support.
\end{proof}

\subsection{An interpretation of the results using connections}
\label{sub:interpretation}
In order to give a geometrical, and perhaps more intuitive, idea of the results
presented, we will describe the theory of Subsection~\ref{sub:hyper} in terms of
a torsion-free linear connection on $TM$. 

Let us consider the setup described at the beginning of Subsection~\ref{sub:hyper};
in addition, we will consider an arbitrary torsion--free (i.e., symmetric) connection
$\nabla$ on the tangent bundle $TM$. We consider the {\em curvature tensor\/} $R$ of $\nabla$
chosen with the following sign convention:
$R(X,Y)=\nabla_X\nabla_Y-\nabla_Y\nabla_X-\nabla_{[X,Y]}$;
the connection $\nabla$ induces a connection $\nabla^*$ in $TM^*$.
For each $p\in TM^*$, the connection $\nabla^*$ produces a direct
sum decomposition of $T_pTM^*$ into a {\em horizontal\/} and a {\em vertical\/}
space in a standard way. The vertical space is identified
with the fiber $T_{\pi(p)}M^*$, and, using $\mathrm d\pi(p)$, the horizontal
space of $T_pTM^*$ is identified with the tangent space $T_{\pi(p)}M$.
These identifications give isomorphisms:
\begin{equation}\label{eq:decompositions}
TTM^*\simeq\pi^*TM\oplus\pi^*TM^*,\quad
T^*TM^*\simeq\pi^*TM^*\oplus\pi^*TM,
\end{equation} 
where $\pi^*$ denotes the pull-back of fiber bundles.
Hence, we have connections $\pi^*\nabla\oplus\pi^*\nabla^*$ and $\pi^*\nabla^*\oplus\pi^*\nabla$ on
$TTM^*$ and
$T^*TM^*$ respectively induced by the above isomorphisms.
Now, given the Hamiltonian function $H$, for each $t\in\R$ we have
a smooth function $H_t$ in (an open subset of) $TM^*$, hence
$\mathrm dH_t$ is a smooth section of $T^*TM^*$. Using the decomposition
\eqref{eq:decompositions}, we write $\mathrm dH=\left(\frac{\partial H}{\partial q},\frac{\partial
H}{\partial p}\right)$, where $\frac{\partial H}{\partial q}$ and $\frac{\partial H}{\partial p}$
denote the restrictions of $\mathrm dH$ to the horizontal and the vertical subspaces
of $TTM^*$ respectively. More in general, we will denote by
$\ddq$ and $\ddp$ respectively the restrictions to the horizontal and
the vertical subspaces of $TTM^*$ of covariant derivatives. One should observe that the connection
in $TTM^*$ is not torsion-free; indeed its torsion is easily computed as:
\[T_p(\rho_1,\rho_2)=(0,R(\mathrm d\pi(\rho_2),\mathrm d\pi(\rho_1))^*p),\]
for $p\in TM^*$ and $\rho_1,\rho_2\in T_pTM^*$.

Using the decomposition \eqref{eq:decompositions} of $TTM^*$, from the Cartan's
identity for exterior differentiation and the symmetry of the connection in $TM$, we
have the following identity for the canonical symplectic form of $TM^*$:
\begin{equation}\label{eq:nuovaomega}
\omega(\rho_1,\rho_2)=[\rho_2]\vertical(\mathrm
d\pi(\rho_1))- [\rho_1]\vertical(\mathrm d\pi(\rho_2)),
\end{equation}
for $p\in TM*$, $\rho_1,\rho_2\in T_pTM^*$, where $[\cdot]\vertical$ denotes projection onto the
vertical space.
Recalling that $\Gamma$ is an integral curve of $\Hf$ in $TM^*$ and that
$\gamma=\pi\circ\Gamma$ is its projection in $M$, the Hamilton  equations
are written as:
\begin{equation}\label{eq:tobelinearized}
\ddt\gamma(t)=\frac{\partial H_t}{\partial
p}(\Gamma(t));\quad
\Ddt\Gamma(t)=-\frac{\partial H_t}{\partial q}(\Gamma(t)),
\end{equation}
where $\frac{\mathrm d}{\mathrm dt}$ denotes standard derivative and $\frac{\mathrm D}{\mathrm
dt}$ denotes covariant derivative. Considering a smooth variation $\Gamma(t,s)$ of $\Gamma$,
$s\in\,]-\varepsilon,\varepsilon\,[$, by integral curves of $\Hf$, with $\Gamma(t,0)=\Gamma(t)$
and with variational vector field $\rho(t)=\dds\Gamma(t,0)$, from
\eqref{eq:tobelinearized} we compute the linearized Hamilton  equations:
\begin{equation}\label{eq:linearized}
 \Ddt v=\frac{\partial^2H_t}{\partial
q\partial p}v+\frac{\partial^2H_t}{\partial
p^2}\alpha;\quad \Ddt\alpha=-\frac{\partial^2H_t}{\partial
q^2}v-R\left(\ddt\gamma,v\right)^*\Gamma-\frac{\partial^2H_t}{\partial
p\,\partial q}\alpha, 
\end{equation}
where $\rho=(v,\alpha)$ is a smooth vector field along $\Gamma$.
We now make an appropriate choice of a symplectic  $\mathfrak L$-trivialization
of $TTM^*$ along $\Gamma$; we consider a parallel trivialization $\{v_1,\ldots,v_n\}$
of $TM$ along $\gamma$, and let $\{\alpha_1,\ldots,\alpha_n\}$ be the corresponding
dual trivialization of $TM^*$ along $\gamma$, which is also parallel.
Using the decomposition of $TTM^*$ in \eqref{eq:decompositions} and formula
\eqref{eq:nuovaomega}, it is easily seen that $\{v_1,\ldots,v_n,\alpha_1,\ldots,\alpha_n\}$
gives a symplectic $\mathfrak L$-trivialization $\phi$ of $TTM^*$ along $\Gamma$.
The coefficients of the symplectic differential system associated to the
Hamiltonian problem by the trivialization $\phi$ are the coordinates in
the basis $\{v_1,\ldots,v_n\}$ of the coefficients of the system 
\eqref{eq:linearized}. Identifying tensors in $TM^*$ with their matrices
in the basis $\{v_1.\ldots,v_n,\alpha_1,\ldots,\alpha_n\}$, we have:
\begin{equation}\label{eq:ABCH}
A=\frac{\partial^2H_t}{\partial q\partial p}, \quad
B=\frac{\partial^2H_t}{\partial p^2},\quad C=-\frac{\partial^2H_t}{\partial
q^2} -\Gamma\circ  R\left(\ddt\gamma,\cdot\right).
\end{equation}

Using this machinery, we are now able to give a geometrical description of the
space $\mathcal K_{\mathcal D}$ introduced in \eqref{eq:KD}. Let's assume that
$\mathcal D$ is a smooth distribution in $M$ defined in a neighborhood of
the image of the curve $\gamma$, such that $\mathcal D_{\gamma(t)}=\mathcal D_t$ 
for all $t$. We say that a curve $\vartheta:[a,b]\to M$ is a {\em solution
of the Hamilton  equations in the direction $\mathcal D$\/} if there exists
a curve $\Theta:[a,b]\to TM^*$, with $\pi\circ\Theta=\vartheta$, such that
the following system is satisfied:
\begin{equation}\label{eq:solD}
\ddt\vartheta(t)=\frac{\partial
H_t}{\partial p}(\Theta(t));\quad
\Ddt\Theta(t)\big\vert_{\mathcal D(\vartheta(t))}=-\frac{\partial H_t}{\partial q}(\Theta(t))
\big\vert_{\mathcal D(\vartheta(t))}.
\end{equation}
An easy computation gives the proof of the following:
\begin{prop}\label{thm:descrKD}
A vector field $\mathfrak v\in
\mathcal H_{\Gamma}$ of Sobolev class $H^2$ is in $\mathcal K_{\mathcal D}$ if and
only if it is the variational vector field of a variation of $\gamma$ by
solutions of the Hamilton  equations in the direction of $\mathcal D$
connecting $P$ and $\gamma(b)$.\qed
\end{prop}

\subsection{The Morse Index Theorem in semi-Riemannian Geometry}
\label{sub:semiRiemannian}
We will now use the results of the previous sections to study the case
of a {\em geodesic Hamiltonian\/} on a semi-Riemannian manifold $(M,g)$,
i.e., $H:TM^*\to\R$ is given by $H(p)=\frac12 g^{-1}(p,p)$ where
$g^{-1}$ denotes the inner product on $TM^*$ induced by $g$. 
Recall that a semi-Riemannian metric on a manifold $M$ is a smooth symmetric
$(2,0)$-tensor field $g$ which is nondegenerate at every point.
It is easy to 
see that $H$ is hyper-regular,  the corresponding Lagrangian
$L:TM\to\R$ is given by $L(v)=\frac12g(v,v)$, and the action functional
is given by:
\[ \vartheta\longmapsto \frac12\int_a^bg(\vartheta',\vartheta')\;\mathrm dt.\]
Let $\nabla$ be the covariant derivative of the Levi--Civita
connection of $g$; in the notation of Subsection~\ref{sub:interpretation}, we
have:
\[\frac{\partial H}{\partial q}=0,\quad\frac{\partial H}{\partial p}=g^{-1}(p),\quad
\frac{\partial^2 H}{\partial q^2}=0,\quad \frac{\partial^2 H}{\partial q\partial p}=0,
\quad \frac{\partial^2 H}{\partial p\partial q}=0,\quad  \frac{\partial^2
H}{\partial p^2}=g^{-1}.\]
By \eqref{eq:tobelinearized}, a solution of the Hamilton  equations
is a curve $\Gamma$ of the form $\Gamma=g(\gamma')$, where $\gamma$
is a geodesic in $(M,g)$. 
As in Subsection~\ref{sub:interpretation}, using a parallel trivialization
of $TM$ along $\gamma$, we obtain an associated symplectic differential
problem whose coefficients \eqref{eq:ABCH} are given by:
\begin{equation}\label{eq:ABCg}
A=0,\quad B=g^{-1},\quad C=g(R(\gamma',\cdot)\,\gamma',\cdot).
\end{equation}
Let $P$ be a submanifold of $M$; the condition
that $\Gamma(a)\in TP^o$ means that the underlying geodesic $\gamma$
starts orthogonally to $P$.
Assumption~\ref{thm:nondegH} in this context means that $g$ is nondegenerate
at $T_{\gamma(a)}P$; by \eqref{eq:linearized}, a solution $v$ of the
linearized Hamilton  equations is a {\em Jacobi field\/}
along $\gamma$, i.e.,
\[v''=R(\gamma',v)\,\gamma'\quad\text{and}\quad \alpha_v=g(v',\cdot).\]
The $\mathcal P$-solutions are the $P$-Jacobi fields along $\gamma$, which are the
Jacobi fields $v$ satisfying the initial conditions:
\[v(a)\in T_{\gamma(a)}P,\quad v'(a)+S^P_{\gamma'(a)}(v(a))\in T_{\gamma(a)}P^\perp,\]
where $S^P_{\gamma'(a)}$ is the second fundamental form of $P$ at $\gamma(a)$ in
the normal direction $\gamma'(a)$. 

A $P$-focal point along $\gamma$ is a point $\gamma(t_0)$, $t_0\in\,]a,b]$, 
such that there exists a non zero $P$-Jacobi field along $\gamma$ vanishing
at $t_0$. For each $t\in\,]a,b]$, we set:
\[\mathcal J[t]=\big\{v(t):v\ \ \text{is a $P$-Jacobi field along}\ \gamma\big\};\]
it is an easy observation that the point $\gamma(t_0)$ is $P$-focal
if and only if $\mathcal J[t]\ne T_{\gamma(t_0)}M$. The {\em multiplicity\/}
of the $P$-focal point $\gamma(t_0)$ is the codimension of $\mathcal J[t_0]$;
$\gamma(t_0)$ is nondegenerate if $g$ is nondegenerate on $\mathcal J[t_0]$ and
the signature of $\gamma(t_0)$ is the signature of the restriction of
$g$ to $\mathcal J[t_0]^\perp$, where $\perp$ denotes the orthogonal complement
with respect to $g$.

If $\gamma(b)$ is not $P$-focal,
the Maslov index $\maslov(\gamma)$ is defined as the Maslov index
$\maslov(\Gamma)$ of the solution of the Hamiltonian $H$; if all the $P$-focal points along
$\gamma$ are nondegenerate, then, by Theorem~\ref{thm:maslov=focal}, $\maslov(\gamma)$ is
equal to the sum of the signatures of all the $P$-focal points along $\gamma$.

The index form $I_\Gamma$ is given by:
\begin{equation}\label{eq:IM}
I_\Gamma(v,w)=\int_a^b\Big[g(v',w')+g(R(\gamma',v)\,\gamma',w)\Big]\;\mathrm
dt-g(S^P_{\gamma'(a)}v(a),w(a)).
\end{equation} 
Let $k=n_-(g)=n_-(\frac{\partial^2H}{\partial
p^2})=n_-(\HL)$; in order to apply Theorem~\ref{thm:main} we need to determine a
$k$-dimensional distribution along $\gamma$ consisting of subspaces where the metric tensor
$g$ is negative definite. Let $\mathcal D$ be one such distribution, generated by
the vector fields $Y_1,\ldots,Y_k$ along $\gamma$. The space $\mathcal K_{\mathcal D}$
defined in \eqref{eq:KD} can be described as:
\begin{equation}\label{eq:newnewK}
\begin{split}
\mathcal K_{\mathcal D}=\big\{v\in\mathcal H_{\Gamma}:&\;g(v',Y_i)\in H^1([a,b];\R)\
\text{and} \\& g(v',Y_i)' =g(v',Y_i')+g(R(\gamma',v)\,\gamma',Y_i),\ \
\forall\,i=1,\ldots,k\big\}.
\end{split}
\end{equation}
The matrices $\mathcal B$, $\mathcal C$ and $\mathcal I$ appearing in
the reduced symplectic system \eqref{eq:symplassoc} are given by:
\begin{equation}\label{eq:newBCI}
\mathcal B_{ij}=g(Y_i,Y_j),\quad\mathcal C_{ij}=g(Y_j',Y_i),\quad
\mathcal I_{ij}=g(Y_i',Y_j')+g(R(\gamma',Y_i)\,\gamma',Y_j).
\end{equation}
We can now restate Theorem~\ref{thm:main} in the form of a Morse Index Theorem
for semi-Riemannian geometry:
\begin{teo}\label{thm:mainsemi}
Let $(M,g)$ be a semi-Riemannian manifold with $n_-(g)=k$, 
let $P\subset M$ be a smooth submanifold
and $\gamma:[a,b]\to M$ be a geodesic with $\gamma'(a)\in TP^\perp$. 
Suppose that $g$ is nondegenerate on $T_{\gamma(a)}P$ and that $\gamma(b)$
is not $P$-focal. Let $Y_1,\ldots,Y_k$ be a family of smooth vector fields
along $\gamma$ generating a $k$-dimensional distribution $\mathcal D$
along $\gamma$ on which $g$ is negative definite. If the reduced symplectic
system \eqref{eq:symplassoc} defined by the matrices \eqref{eq:newBCI}
has no focal instants, then the index of $I_\Gamma$ (\eqref{eq:IM})
on the space $\mathcal K_{\mathcal D}$ defined in \eqref{eq:newnewK}
is given by:
\[\hbox to\hsize{\hfill$\displaystyle\phantom{\qed} n_-(I_\Gamma\big\vert_{\mathcal
K_{\mathcal D}})=\maslov(\gamma)+n_-(g
\big\vert_{T_{\gamma(a)}P}).\hfill\qed$}\]
\end{teo}
Observe that if $(M,g)$ is Riemannian, then Theorem~\ref{thm:mainsemi}
is the Morse Index Theorem for Riemannian geodesics with variable initial
points; if $(M,g)$ is Lorentzian and $\gamma$ is a timelike geodesic, then
it is not hard to see that, taking $Y_1=\gamma'$ in Theorem~\ref{thm:mainsemi},
we obtain the Timelike Morse Index Theorem in Lorentzian geometry
(see \cite{BEE, PT}). Theorem~\ref{thm:mainsemi} generalizes \cite[Theorem~6.1]{GMPT}.

We present some examples where the Theorem~\ref{thm:mainsemi} applies.
\begin{example}\label{thm:exaJacobi}
If $Y_1,\ldots,Y_k$ are Jacobi fields such that $g(Y_i',Y_j)$ is symmetric,
then the reduced symplectic system \eqref{eq:symplassoc} has no focal
instants (see Example~\ref{thm:example1}). In this case, the space
$\mathcal K_{\mathcal D}$ can be easily described as:
\[\mathcal K_{\mathcal D}=\big\{v\in\mathcal H_{\Gamma}:
g(v',Y_i)-g(v,Y_i')=\text{constant},\ \ i=1,\ldots,k\big\}.\]
\end{example}

\begin{example}\label{thm:exaKilling}
Let $G$ be a $k$-dimensional Lie group acting on $M$ by isometries
with no fixed points, or more in general, having only discrete
isotropy groups. Suppose that $g$ is negative definite on the orbits of $G$.
If $\gamma'(a)$ is orthogonal to the orbit of the commutator subgroup
$[G,G]$ (for instance if $G$ is abelian), then Theorem~\ref{thm:mainsemi}
can be applied by considering the distribution $\mathcal D$ tangent
to the orbits of $G$. Observe indeed that $\mathcal D$ is generated
by $k$ linearly independent Killing vector fields $Y_1,\ldots,Y_k$
on $M$, which therefore restrict to Jacobi fields along any geodesic.
Then, one falls into the situation of Example~\ref{thm:exaJacobi}
by observing  that the symmetry of $g(Y_i',Y_j)$ follows from
the orthogonality with the orbits of $[G,G]$:
\[g(Y_i',Y_j)-g(Y_i,Y_j')=-g(\nabla_{Y_j}Y_i,\gamma')+g(\nabla_{Y_i}Y_j,\gamma')=
g([Y_i,Y_j],\gamma')=0.\]
In this situation, the space $\mathcal K_{\mathcal D}$ can be described
as the space of variational vector fields along $\gamma$ corresponding to
variations of $\gamma$ by curves that are {\em geodesics along $\mathcal D$}, i.e.,
whose second derivatives are orthogonal to $\mathcal D$ 
(see Proposition~\ref{thm:descrKD}).
\end{example}

\begin{example}\label{thm:exacurvature}
Suppose that the bilinear form $g(R(\gamma',\cdot)\,\gamma',\cdot)$ is
negative semidefinite along the geodesic $\gamma$. Then, Theorem~\ref{thm:mainsemi}
can be applied by considering any $k$-dimen\-sional parallel distribution $\mathcal D$ along
$\gamma$ where $g$ is negative definite (see Example~\ref{thm:example2}).
If $k=1$, i.e., if $(M,g)$ is Lorentzian, this  is a condition on the sign
of the sectional curvature of timelike planes containing $\gamma'$.
\end{example}
We conclude with the observation that from Theorem~\ref{thm:indexthmPQ}
one can easily obtain a semi-Riemannian Morse Index Theorem for
geodesics starting and ending orthogonally to two submanifolds of $M$.
 \end{section}




\begin{thebibliography}{99}

\bibitem{AM} R.\ Abraham,  J.\ E.\ Marsden, {\em Foundations of Mechanics}, 2nd
Edition, Benjamin/Cummings, Advanced Book Program, Reading, Mass., 1978.

\bibitem{Ar} V.\ I.\ Arnol'd, {\em Characteristic Class Entering in Quantization
Conditions}, Funct.\ Anal.\ Appl.\ {\bf1} (1967), 1--13.

\bibitem{BEE} J.\ K.\ Beem, P.\ E.\ Ehrlich, K.\ L.\ Easley,
{\em Global Lorentzian Geometry},
Marcel Dekker, Inc., New York and Basel, 1996.

\bibitem{Bott} R.\ Bott, {\em On the Iteration of Closed Geodesics and
the Sturm Intersection Theory}, Commun.\ Pure Appl.\ Math.\ {\bf9} (1956), 171--206.

\bibitem{CL} E.\ A.\ Coddington, N.\ Levinson, {\em Theory of Ordinary Differential
Equations}, McGraw--Hill Book Company, New York, Toronto, London, 1955.

\bibitem{dC} M.\ do Carmo, {\em Riemannian Geometry}, Birkh\"auser, Boston, 1992.

\bibitem{CZ} C.\ Conley, E.\ Zehnder, {\em Morse-type Index Theory for Flows
and Periodic Solutions for Hamiltonian Equations}, Comm.\ Pure Appl.\ Math.\ {\bf37} (1984),
207--253.

\bibitem{Duis} J.\ J.\ Duistermaat, {\em On the Morse Index in Variational Calculus},
Adv.\ in Math.\ {\bf 21} (1976), 173--195.

\bibitem{Ed} H.\ M.\ Edwards, {\em A Generalized Sturm Theorem}, Ann.\ of Math.\ {\bf80}
(1964), 22--57.

\bibitem{Ek} I.\ Ekeland, {\em An Index Theory for Periodic Solutions of Convex Hamiltonian
Systems}, Proc.\ of Symposia in Pure Math.\ vol.\ 45 (1986), Part I, 395--423.

\bibitem{GMPT} F.\ Giannoni, A.\ Masiello, P.\ Piccione,
D.\ Tausk, {\em A Generalized Index Theorem for Morse--Sturm Systems
and Applications to semi-Riemannian Geometry}, preprint 1999 ({\tt LANL math.DG/9908056}),
to appear in the Asian Journal of Mathematics.

\bibitem{GS} V.\ Guillemin, S.\ Sternberg, {\em Geometric Asymptotics}, Mathematical Surveys
and Monographs n.\ 14, AMS, Providence RI, 1990.

\bibitem{Hel1} A.\ D.\ Helfer, {\em Conjugate Points on Spacelike Geodesics
or Pseudo-Self-Adjoint Morse-Sturm-Liouville Systems}, Pacific J.\ Math.\ {\bf164}, n.\ 2
(1994), 321--340.


\bibitem{Hirsch} M.\ W.\ Hirsch, {\em Differential Topology}, Springer--Verlag, 1976.

\bibitem{Kal} D.\ Kalish, {\em The Morse Index Theorem where the Ends are Submanifolds},
Trans.\ Am.\ Math.\ Soc.\ {\bf308}, n.\ 1 (1988), 341--348.

\bibitem{Long} Y.\ Long, {\em A Maslov-type Index Theory for Symplectic Paths},
Top.\ Meth.\ Nonlin.\ Anal.\ {\bf 10} (1997), 47--78.

\bibitem{Ma} A.\ Masiello, {\em
Variational Methods in Lorentzian Geometry}, 
Pitman Research Notes in Mathematics {\bf 309}, 
Longman, London 1994.

\bibitem{MW} J.\ Mawhin, M.\ Willem, {\em Critical Point Theory and
Hamiltonian Systems}, Springer--Verlag, Berlin, 1989.

\bibitem{MPT} F.\ Mercuri, P.\ Piccione, D.\ Tausk,
{\em Stability of the Focal and the Geometric Index in semi-Riemannian
Geometry via the Maslov Index}, Technical Report RT-MAT 99-08, 
Mathematics Department, University of S\~ao Paulo, Brazil, 1999. ({\tt LANL math.DG/9905096})

\bibitem{M} J.\ Milnor, {\em Morse Theory}, Princeton Univ.\ Press, Princeton, 1969.


\bibitem{Perlick1999} Perlick, V.
        {\em Ray optics, Fermat's principle and applications
        to general relativity}, Sprin\-ger Lecture Notes in
        Physics, m61.


\bibitem{PPT} V.\ Perlick, P.\ Piccione, D.\ Tausk,
{\em A Variational Principle for Light Rays in
a General  Relativistic Plasma}, preprint 1999.

\bibitem{PT} P.\ Piccione, D.\ V.\ Tausk,  {\em A Note on the Morse Index Theorem for Geodesics
between Submanifolds in semi-Riemannian Geometry}, J.\  Math.\  Phys.\ {\bf40}, vol.\ 12
(1999),  6682--6688.

\bibitem{CRASP} P.\ Piccione, D.\ V.\ Tausk, {\em The Maslov Index and a Generalized Morse
Index Theorem for Non Positive Definite Metrics}, Comptes Rendus de l'Acad\'emie
de Sciences de Paris, vol.\  331, 5 (2000),
385--389. 



\bibitem{SZ} D.\ Salamon, E.\ Zehnder, {\em Morse Theory for Periodic
Solutions of Hamiltonian Systems and the Maslov Index}, Commun.\ Pure Appl.\ Math.\ 
{\bf45} (1992), 1303--1360.


\bibitem{Tre} F.\ Treves, {\em Introduction to Pseudodifferential Operators}, Plenum,
New York, 1982.

\end{thebibliography}
\end{document}